\documentclass[12pt,leqno]{article}
\def\reals{\hbox{\sl I\kern-.18em R \kern-.3em}}
\def\quats{\hbox{\sl I\kern-.18em H \kern-.3em}}
\def\ints{\hbox{\sl Z\kern-.4em Z \kern-.3em}}
\def\nats{\hbox{\sl I\kern-.16em N \kern.05em}}
\def\rats{\hbox{\sl Q \kern-.83em\vrule height.59em depth0em \kern.87em}}
\def\complexes{\hbox{\sl\kern.50em I\kern-.50em C \kern.05em}}
\newtheorem{thm}{Theorem}[section]
\newtheorem{prop}[thm]{Proposition}
\newtheorem{coro}[thm]{Corollary}
\newtheorem{lemm}[thm]{Lemma}
\newtheorem{remark}[thm]{Remark}
\def\proof{\noindent {\bf Proof: }}
\def\np{\noindent}
\def\endpf{$\|$  \medskip}
\newcommand{\ov}[1]{#1^\ast}
\newcommand{\pos}{B_n^+}

\def\cy{\mbox{\bf c}}
\def\dy{\mbox{\bf d}}
\def\inf{{\rm inf}}
\def\sup{{\rm sup}}

\begin{document}

\title{A new approach to the word and conjugacy problems in the braid
groups}
\author{Joan S. Birman\footnote{Partially supported by NSF Grant
94-02988. This paper was completed during a visit by the first author
to MSRI. She thanks MSRI for its hospitality and partial support, and
thanks Barnard
College for its support under a Senior Faculty Research Leave.},  Ki Hyoung
Ko\footnote{This work was initiated during the second author's sabbatical
visit to Columbia University in 1995-6. He thanks Columbia University for
its hospitality during that visit.} and Sang Jin Lee.}

\maketitle
\centerline {11/97 revision of 5/97 draft; to appear in
ADVANCES IN MATHEMATICS}
\centerline{previous title:``A new approach to the word problem in the braid groups}
\begin{abstract}
\noindent A new presentation of the $n$-string braid group $B_n$ is studied.
Using it, a new solution to the word problem in $B_n$
is obtained which retains most of the  desirable features of the
Garside-Thurston solution, and at the same time makes possible
certain computational improvements. We also give a related solution to
the conjugacy problem, but the improvements
in its complexity are not clear at this writing.
\end{abstract}
\section{Introduction}

In the foundational manuscript \cite{Artin} Emil Artin introduced the
sequence
of braid groups $B_n, \ n=1,2,3,\dots$ and proved that $B_n$ has a
presentation
with $n-1$ generators
$\sigma_1,\sigma_2,\ldots,\sigma_{n-1}$ and defining relations:
\begin{equation}
\label{artinrelation1} \sigma_t \sigma_s =  \sigma_s \sigma_t \ \ \ \ {\rm
if
} \ \ \ |t-s|>1.
\end{equation}
\begin{equation}
\label{artinrelation2}
\sigma_t\sigma_s\sigma_t = \sigma_s\sigma_t\sigma_s \ \ \ \ {\rm if
} \ \ \ |t-s| = 1.
\end{equation}
The word problem in $B_n$ was posed by Artin in \cite{Artin}.
His solution was based on his knowledge
of the structure of the kernel of the map $\phi$ from $B_n$ to the
symmetric group $\Sigma_n$ which sends the generator $\sigma_i$ to the
transposition $(i,i+1)$.  He used the group-theoretic properties of
the kernel of $\phi$ to
put a braid into a normal form called a `combed braid'. While nobody has
investigated the matter, it seems intuitively clear that Artin's
solution is exponential in the length of a word in the generators $\sigma_1,
\dots,\sigma_{n-1}$.

The conjugacy problem in $B_n$ was also posed
in \cite{Artin}, also its importance for the problem of
recognizing knots and links algorithmically was noted, however it took 43
years before progress was made.  In a different, but equally foundational
manuscript  \cite{Garside}  F. Garside
discovered a new solution to the word problem (very different from
Artin's) which then led him to a related solution to the conjugacy
problem. In Garside's solution one focusses not on the kernel of $\phi$, but
on its image, the symmetric group $\Sigma_n$.
Garside's solutions to both the word and conjugacy
problem are exponential in both word length and braid index.

The question of the speed of Garside's algorithm for the word problem
was first raised by Thurston. His contributions, updated to reflect
improvements obtained after his widely circulated preprint appeared, are
presented in Chapter 9 of
\cite{Epstein}.  In \cite{Epstein} Garside's algorithm is modified by
introducing new ideas, based upon the fact that braid groups are
biautomatic, also that $B_n$ has a partial ordering which gives it the
structure of a lattice.  Using these facts it is proved in  \cite{Epstein}
that there exists  an algorithmic solution to the word problem which is
${\cal O}(|W|^2 n \ log \ n)$, where
$|W|$ is word length. See, in particular, Proposition 9.5.1 of
\cite{Epstein}, our discussion at the beginning of $\S 4$ below, and
Remark \ref{complexity} in $\S 4$.  While the same general set of ideas
apply equally well to the conjugacy problem
\cite{Elrifai:Morton}, similar sharp estimates  of complexity have not been
found because the combinatorial complications present a new level of
difficulty.

A somewhat different question is the {\it shortest word problem}, to find
a representative of the word class which has shortest length in
the Artin generators. It was proved in
\cite{PR} that this problem in $B_n$
is at least as hard as  an NP-complete problem. Thus, if one could
find a polynomial time algorithm to solve the shortest word problem one
would have proved that P=NP.

Our contribution to this set of ideas is to introduce a new and very
natural set of generators for $B_n$ which includes the Artin generators as
a subset. Using the new generators we will be able to solve the word
problem in much the same way as Garside and Thurston solved it,
moreover our solution generalizes to a related solution to the conjugacy
problem which is in the spirit of that of \cite{Elrifai:Morton}.
The detailed combinatorics in our work are, however, rather
different from those in \cite{Elrifai:Morton} and \cite{Epstein}. Our
algorithm solves the word problem in ${\cal O}(|W|^2n)$. Savings in actual
running time (rather than complexity) also occur, because a word written in
our generators is generally shorter by a factor of $n$ than a word in the
standard generators which represents the same element (each
generator $a_{ts}$ in our work replaces a word of length $2(t-s)-1$, where
$n>t-s>0$ in the Artin generators),  also the positive
part is shorter by a factor of $n$ because the new generators lead to a
new and shorter `fundamental word' $\delta$ which replaces Garside's
famous $\Delta$.

Our solution to both the word and conjugacy problems generalizes the
work  of Xu
\cite{Xu} and of Kang, Ko and Lee \cite{KKL}, who succeeded in finding
polynomial time algorithms for the word and conjugacy problems and also
for the shortest word problem in $B_n$ for $n=3$ and
$4$.     The general case appears to be more subtle than the cases $n=3$
and $4$, however
polynomial time solutions to the three problems for every $n$ do not seem to
be totally out of reach, using our generators.

In the three references \cite{Elrifai:Morton}, \cite{Epstein} and
\cite{Garside} a central role is played by {\em
positive braids}, i.e. braids which are positive powers of the generators.
Garside introduced  the {\em fundamental
braid} $\Delta$:
\begin{equation}
\Delta = (\sigma_1\sigma_2\dots\sigma_{n-1})
(\sigma_1\sigma_2\dots\sigma_{n-2})\dots(\sigma_1\sigma_2)(\sigma_1).
\end{equation}
He showed that every element ${\cal W} \in B_n$ can be represented
algorithmically by a word $W$ of the form $\Delta^rP$, where $r$ is an
integer and $P$ is a positive word, and $r$ is maximal for all such
representations. However his  $P$ is non-unique up to a finite set of
equivalent words which represent the same element ${\cal P}$. These can
all be found algorithmically, but the list is very long. Thus instead of a
unique normal form one has a fixed
$r$ and a finite set of positive words which represent ${\cal P}$.
Thurston's improvement was to show that ${\cal P}$ can in fact be factorized
as a  product $P_1P_2\dots P_s$, where each
$P_j$ is a special type of
positive braid which is known as `permutation braid'. Permutation braids
are determined uniquely by their associated permutations, and Thurston's
normal form is a unique representation of this type in which the integer
$s$ is minimal for all representations of ${\cal P}$ as a product of
permutation braids. Also, in each subsequence $P_iP_{i+1}\dots P_s,
i=1,2,\dots,s-1,$ the permutation braid $P_i$ is the longest possible
permutation braid in a factorization of this type. The subsequent work of
Elrifai and Morton \cite{Elrifai:Morton} showed that there is a related
algorithm which simultaneously maximizes $r$ and minimizes $s$ within each
conjugacy class. The set of all products
$P_1P_2\dots P_s$ which do that job (the {\em super summit set}) is
finite, but it is not well understood.

Like Artin's, our generators are
braids in which exactly one pair of strands crosses, however the images of
our  generators in
$\Sigma_n$ are {\it arbitrary} transpositions $(i,j)$ instead of simply
adjacent transpositions $(i,i+1)$.  For each $t,s$ with $n\geq
t>s\geq 1$ we consider the element of $B_n$ which is defined by:
\begin{equation}
\label{newgenerators}
a_{ts}=(\sigma_{t-1}\sigma_{t-2}\cdots\sigma_{s+1})\sigma_s
(\sigma_{s+1}^{-1}\cdots\sigma_{t-2}^{-1}\sigma_{t-1}^{-1}).
\end{equation}
so that our generators include the Artin generators (as a proper subset
for $n\geq 3$). The braid $a_{ts}$ is depicted in Figure 1(a).
Notice that $a_{21},a_{32}\dots$ coincide with $\sigma_1,\sigma_2,\dots$

\np The braid $a_{ts}$
is an elementary interchange of the $t^{th}$ and $s^{th}$ strands, with all
other
strands held fixed, and with the convention that the strands being
interchanged
pass in front of all intervening strands. We call them {\it band
generators} because
they suggest a disc-band decomposition of a surface bounded by a closed
braid.

We introduce a new {\it fundamental word}:
\begin{equation}
\label{eqn:delta}
 \delta  = a_{n(n-1)}a_{(n-1)(n-2)}\cdots a_{21} =
\sigma_{n-1}\sigma_{n-2}\dots\sigma_2\sigma_1.
\end{equation}
The reader who is familiar with the mathematics of braids will recognize
that $\Delta^2 = \delta^n$ generates the center of $B_n$. Thus $\Delta$
may be thought
of as the `square root' of the center, whereas $\delta$ is the `nth root'
of the
center. We will prove that each element ${\cal W}\in B_n$ may be
represented (in terms of the band generators) by a unique word $W$ of the
form:
\begin{equation}
W= \delta^jA_1A_2\cdots A_k,
\end{equation}
where $A = A_1A_2\cdots A_k$ is positive, also
$j$ is maximal and $k$ is minimal for
all such representations, also the $A_i$'s are positive braids
which are determined uniquely by their associated permutations. We will
refer to Thurston's braids $P_i$ as {\em permutation braids}, and to our
braids $A_i$ as  {\em canonical factors}.

\begin{picture}(420,208)(60,585)
\put(100,585){\makebox(0,0)[cb]{\smash{\normalsize\rm (a)}}}
\put(95,610){\makebox(0,0)[lb]{\smash{\Large\rm.}}}
\put(95,615){\makebox(0,0)[lb]{\smash{\Large\rm.}}}
\put(95,620){\makebox(0,0)[lb]{\smash{\Large\rm.}}}
\put(95,755){\makebox(0,0)[lb]{\smash{\Large\rm.}}}
\put(95,760){\makebox(0,0)[lb]{\smash{\Large\rm.}}}
\put(95,765){\makebox(0,0)[lb]{\smash{\Large\rm.}}}
\put(80,675){\makebox(0,0)[lb]{\smash{\Large\rm.}}}
\put(80,680){\makebox(0,0)[lb]{\smash{\Large\rm.}}}
\put(80,685){\makebox(0,0)[lb]{\smash{\Large\rm.}}}
\put(80,690){\makebox(0,0)[lb]{\smash{\Large\rm.}}}
\put(80,695){\makebox(0,0)[lb]{\smash{\Large\rm.}}}
\put(80,700){\makebox(0,0)[lb]{\smash{\Large\rm.}}}
\put(65,733){\makebox(0,0)[rb]{\smash{\small\rm$t$-th}}}
\put(65,643){\makebox(0,0)[rb]{\smash{\small\rm$s$-th}}}
\thinlines
\put( 60,605){\line( 1, 0){ 80}}
\put( 60,630){\line( 1, 0){ 80}}
\put( 60,750){\line( 1, 0){ 80}}
\put( 60,770){\line( 1, 0){ 80}}
\put( 70,645){\line( 1, 0){ 25}}
\put(95,645){\line( 1, 6){  5.811}}
\put(110,735){\line(-1,-6){  5,811}}
\put(110,735){\line( 1, 0){ 30}}
\put(70,735){\line( 1, 0){ 25}}
\put(95,735){\line( 1,-6){ 15}}
\put(110,645){\line( 1, 0){ 30}}
\put( 60,660){\line( 1, 0){ 32}}
\put(100,660){\line( 1, 0){  5}}
\put(113,660){\line( 1, 0){ 27}}
\put( 60,720){\line( 1, 0){ 32}}
\put(100,720){\line( 1, 0){  5}}
\put(113,720){\line( 1, 0){ 27}}
\put(340,585){\makebox(0,0)[cb]{\smash{\normalsize\rm (b)}}}
\put(195,758){\makebox(0,0)[rb]{\smash{\small\rm$t$-th}}}
\put(195,668){\makebox(0,0)[rb]{\smash{\small\rm$s$-th}}}
\put(195,613){\makebox(0,0)[rb]{\smash{\small\rm$r$-th}}}
\put(200,760){\line( 1, 0){ 15}}
\put(215,760){\line( 1,-6){ 15}}
\put(230,670){\line( 1, 0){ 10}}
\put(240,670){\line( 1,-5){ 11}}
\put(251,615){\line( 1, 0){ 15}}
\put(200,670){\line( 1, 0){ 15}}
\put(215,670){\line( 1, 6){  6}}
\put(230,760){\line(-1,-6){  6}}
\put(230,760){\line( 1, 0){ 36}}
\put(200,615){\line( 1, 0){ 40}}
\put(240,615){\line( 1, 5){  4}}
\put(251,670){\line(-1,-5){  4}}
\put(251,670){\line( 1, 0){ 15}}
\put(281,685){\makebox(0,0)[cb]{\smash{\Large $=$}}}
\put(296,760){\line( 1, 0){ 15}}
\put(311,760){\line( 1,-6){ 24.18}}
\put(335.18,615){\line( 1, 0){ 50}}
\put(296,670){\line( 1, 0){ 19}}
\put(321,670){\line( 1, 0){  3}}
\put(328,670){\line( 1, 0){ 17}}
\put(345,670){\line( 1, 6){  6}}
\put(360,760){\line(-1,-6){  6}}
\put(360,760){\line( 1, 0){ 15}}
\put(296,615){\line( 1, 0){ 15}}
\put(311,615){\line( 1, 6){ 11}}
\put(335.18,760){\line(-1,-6){11}}
\put(335.18,760){\line( 1, 0){ 9.82}}
\put(345,760){\line( 1,-6){ 15}}
\put(360,670){\line( 1, 0){ 15}}
\put(390,685){\makebox(0,0)[cb]{\smash{\Large $=$}}}
\put(405,760){\line( 1, 0){ 36}}
\put(441,760){\line( 1,-6){ 24.18}}
\put(465.18,615){\line( 1, 0){ 15}}
\put(405,670){\line( 1, 0){ 15}}
\put(420,670){\line( 1,-5){ 11}}
\put(431,615){\line( 1, 0){ 10}}
\put(441,615){\line( 1, 6){ 11}}
\put(465.18,760){\line(-1,-6){ 11}}
\put(465.18,760){\line( 1, 0){ 15}}
\put(405,615){\line( 1, 0){ 15}}
\put(420,615){\line( 1, 5){ 4}}
\put(431,670){\line(-1,-6){ 4}}
\put(431,670){\line( 1, 0){ 17}}
\put(452,670){\line( 1, 0){  3}}
\put(459,670){\line( 1, 0){ 21}}
\end{picture}

\

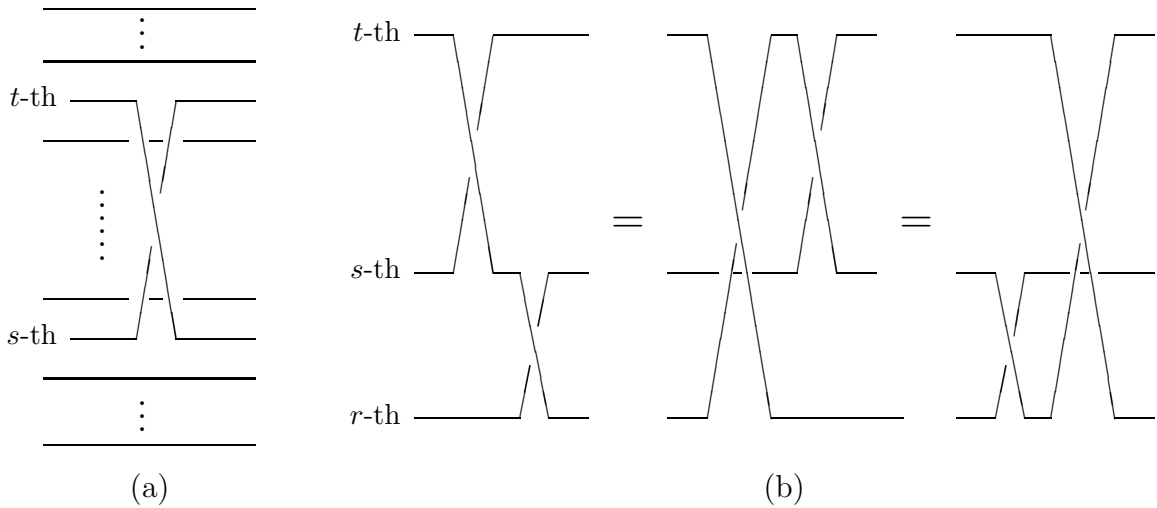
\centerline{\bf Figure 1: The band generators and relations between them}

\

Let ${\cal W}$ be an arbitrary element of $B_n$ and let $W$ be a word in
the band generators which represents it. We are able to analyze the speed of
our algorithm for the word problem, as a function of both the word length
$|W|$  and braid index $n$.  Our main result is a
new algorithmic solution
to the word problem (see $\S$4 below). Its computational complexity,
which is analysed carefully in $\S$4 of this paper,
is an improvement over that given in \cite{Epstein} which is the best
among the known algorithms.
Moreover our work offers certain other advantages, namely:
\begin{enumerate}
\item The number of
distinct permutation braids is $n!$, which grows faster than $k^n$ for
any $k\in\reals^+$.  The number of
distinct canonical factors is the $n^{th}$ Catalan number
${\cal C}_n = (2n)!/n!(n+1)!$, which is bounded above by $4^n$.
The reason for this reduction is the fact that the canonical factors can be
decomposed nicely into parallel, descending cycles (see
Theorem~\ref{thm:subwords}).
The improvement in the complexity of the word algorithm is
a result of the fact that the canonical factors are very simple.
We think that they reveal beautiful new structure in the braid group.
\item Since our generators include
the Artin generators, we may assume in both cases that we begin with a word
$W$ of length $|W|$ in the Artin generators. Garside's $\Delta$ has
length $(n-1)(n-2)/2$, which implies that the
word length $|P|$ of the positive word $P =P_1P_2\cdots P_q$ is roughly
$n^2|W|$. On the other hand,  our $\delta$
has word length $n-1$, which implies that the length $|A|$ of the product
$A = A_1A_2\cdots A_k$ is roughly $n|W|$.
\item Our work, like that in \cite{Epstein}, generalizes to the conjugacy
problem. We conjecture that our solution to that problem is polynomial in
word length, a matter which we have not settled at this writing.
\item Our solution to the word problem suggests a related solution to
the shortest word problem.
\item It has been noted in conversations with A. Ram that our work
ought to generalize to other Artin groups
with finite Coxeter groups. This may be of interest in its own right.
\end{enumerate}

Here is an outline of the paper. In $\S$2 we find a presentation for
$B_n$ in terms of the new generators and show that there is  a natural
semigroup
$B_n^+$  of positive words which is determined by the presentation.
We prove that every element in $B_n$ can be represented in the form
$\delta^tA$, where $A$ is a positive word. We then prove (by a long
computation) that
$B_n^+$ embeds in $B_n$, i.e. two positive words in $B_n$ represent the
same element of $B_n$ if and only if their pullbacks to $B_n^+$ are equal
in $B_n^+$. We note (see Remark 2.8) that our generators and Artin's are
the only ones in a class studied in \cite{Sergiescu} for which such an
embedding theorm holds. In
$\S$3 we use these ideas to find normal forms for words in $B_n^+$, and so
also for words in $B_n$. In $\S$4 we give our algorithmic
solution to the word problem and study its complexity. In $\S$5 we
describe very briefly how our work generalizes to the conjugacy problem.

\

\begin{remark}
\label{remark:literature}  {\rm In the article \cite{Dehornoy}
P. Dehornoy gives an algorithmic solution to the word problem which
is based upon the existence, proved in a different paper by the same
author, of an order structure on $B_n$.  His methods seem quite different
from ours and from those in the other papers we have cited, and not in a
form where precise comparisons are possible. Dehornoy does not discuss the
conjugacy problem, and indeed his methods do not seem to generalize to
the conjugacy problem.}
\end{remark}

\

\noindent{\bf Acknowledgements}  We thank Marta Rampichini for her
careful reading of earlier versions of this manuscript, and her thoughtful
questions. We thank Hessam Hamidi-Tehrani for pointing out to us
the need to clarify our calculations of computational complexity.

\newpage

\section{The semigroup of positive braids}
We begin by finding a presentation for $B_n$ in terms of the new
generators. We
will use the symbol $a_{ts}$ whenever there is no confusion about the two
subscripts,
and symbols such as $a_{(t+2)(s+1)}$ when there might be confusion
distinguishing
between the first and second subscripts. Thus $a_{(t+1)t}=\sigma_t$.

\begin{prop}\label{thm:presentation}
$B_n$ has a presentation with generators $\{a_{ts} ; \ n\geq t>s\geq
1 \}$ and with defining relations
\end{prop}
\begin{equation}
\label{eqn:bandrelation1}
a_{ts}a_{rq}=a_{rq}a_{ts}\quad\mbox{\rm if } \ (t-r)(t-q)(s-r)(s-q)>0
\end{equation}
\begin{equation}
\label{eqn:bandrelation2}
a_{ts}a_{sr}=a_{tr}a_{ts}=a_{sr}a_{tr} \ \ {\rm for \  all} \ t,s,r
 {\rm \  with} \ \ n\geq t>s>r \geq 1.
\end{equation}

\begin{remark}
\label{remark:bandgenerator relations}
{\rm Relation
(\ref{eqn:bandrelation1}) asserts that $a_{ts}$ and $a_{rq}$ commute if
$t$ and $s$ do
not separate $r$ and $q$. Relation (\ref{eqn:bandrelation2}) expresses a
type of
`partial' commutativity in the case when $a_{ts}$ and $a_{rq}$ share a
common strand.
It tells us that if the product $a_{ts}a_{sr}$ occurs in a braid word,
then we may move
$a_{ts}$ to the right (resp. move $a_{sr}$ to the left) at the expense of
increasing
the first subscript of $a_{sr}$ to $t$ (resp. decreasing the second
subscript of
$a_{ts}$ to $r$. )}
\end{remark}

\proof We begin with Artin's presentation for $B_n$, using generators
$\sigma_1,\dots,\sigma_{n-1}$ and relations (\ref{artinrelation1}) and
(\ref{artinrelation2}).  Add the new generators $a_{ts}$ and the relations
(\ref{newgenerators}) which define them in terms of the $\sigma_i$'s.
Since we know
that relations (\ref{eqn:bandrelation1}) and (\ref{eqn:bandrelation2})
are described by
isotopies of braids, depicted in Figure 1(b), they must
be
consequences of (\ref{artinrelation1}) and (\ref{artinrelation2}), so we
may add them
too.

In the special case when $t=s+1$ relation (\ref{newgenerators}) tells us
that
$a_{(i+1)i}=\sigma_i$, so we may omit the generators
$\sigma_1,\dots,\sigma_{n-1}$, to obtain a presentation with generators
$a_{ts}$, as described in the theorem. Defining relations are now
(\ref{eqn:bandrelation1}), (\ref{eqn:bandrelation2}) and:
\begin{equation}
\label{artinrelation3}
a_{(t+1)t}a_{(s+1)s} = a_{(s+1)s}a_{(t+1)t}\quad\mbox{if } |t-s|>1
\end{equation}
\begin{equation}
\label{artinrelation4}
a_{(t+1)t}a_{(s+1)s}a_{(t+1)t} = a_{(s+1)s}a_{(t+1)t}a_{(s+1)s}
{\rm \ \ if \ \ } |t - s| = 1
\end{equation}
\begin{equation}
\label{newgenerators2}
a_{ts}=(a_{t(t-1)}a_{(t-1)(t-2)}\cdots a_{(s+2)(s+1)})a_{(s+1)s}
(a_{t(t-1)}a_{(t-1)(t-2)}\cdots a_{(s+2)(s+1)})^{-1}
\end{equation}
 Our task is to prove that
(\ref{artinrelation3}), (\ref{artinrelation4}) and (\ref{newgenerators2})
are
consequences of (\ref{eqn:bandrelation1}) and
(\ref{eqn:bandrelation2}).

Relation (\ref{artinrelation3}) is nothing more than a special case of
(\ref{eqn:bandrelation1}).
As for (\ref{artinrelation4}), by symmetry we may assume that $t=s+1$.
Use (\ref{eqn:bandrelation2}) to replace $a_{(s+2)(s+1)}a_{(s+1)s}$ by
$a_{(s+1)s}a_{(s+2)s}$, thereby reducing (\ref{artinrelation4}) to
$a_{(s+2)s}a_{(s+2)(s+1)}=a_{(s+2)(s+1)}a_{(s+1,s)}$, which is a
special case of (\ref{eqn:bandrelation2}).
Finally, we consider (\ref{newgenerators2}). If $t=s+1$ this relation is
trivial, so we may assume that $t > s+1$. Apply (\ref{eqn:bandrelation2})
to change
the center pair $a_{(s+2)(s+1)}a_{(s+1)s}$ to $a_{(s+2)s}a_{(s+2)(s+1)}$.
If
$t>s+2$  repeat this move on the new pair $a_{(s+3)(s+2)}a_{(s+2)s}$.
Ultimately, this process will move the original center letter $a_{(s+1)s}$
to the
leftmost position, where it becomes $a_{ts}$. Free cancellation
eliminates everything to its right, and we are done.  \endpf

A key feature which the new presentation shares with the old is that the
relations have all been expressed as relations between positive powers of
the generators, also the relations all preserve word length. Thus our
presentation also determines a presentation for a semigroup. A word in
positive powers of the generators  is called a {\em positive word}.  Two
positive words are said to be {\em positively equivalent}  if one can be
transformed into the other by a sequence of  positive words such that each
word of the sequence is obtained from the preceding one by a single direct
application  of a defining relation in (\ref{eqn:bandrelation1}) or
(\ref{eqn:bandrelation2}). For two positive words $X$ and $Y$, write
$X\doteq Y$
if they  are positively equivalent. Positive words that are positively
equivalent have the same word length since all of defining relations
preserve the word length.  We use the symbol $B_n^+$ for the monoid of
positive
braids, which can be defined by the generators and
relations in Theorem~\ref{thm:presentation}.
Thus $\pos$ is the set of positive words modulo positive
equivalence. Our next goal is to prove that the principal
theorem
of \cite{Garside} generalizes to our new presentation, i.e. that the
monoid of
positive braids embeds in the braid group $B_n$.  See Theorem
\ref{thm:embedding}
below.

Before we can begin  we need to establish  key properties of the fundamental
braid
$\delta$. Let $\tau$ be the inner automorphism of
$B_n$ which is induced by
conjugation by $\delta$.

\begin{lemm}
\label{lemm:delta}
Let $\delta$ be the fundamental braid. Then:
\begin{enumerate}
\item [{\rm (I)}] $\delta = a_{n(n-1)}a_{(n-1)(n-2)}\cdots a_{21}$ is
positively
equivalent to a word that begins or ends with any given generator $a_{ts},
\ \ n
\geq t > s \geq 1$.  The explicit expressions are:
\begin{enumerate}
\item [ ]
$\delta \doteq (a_{ts})(a_{n(n-1)}\cdots
a_{(t+2)(t+1)}a_{(t+1)s}a_{s(s-1)}\cdots a_{21})(a_{t(t-1)}\cdots
a_{(s+2)(s+1)}) $
\item [ ]
$\delta \doteq  (a_{n(n-1)}a_{(n-1)(n-2)}\cdots
a_{(t+1)t}a_{t(s-1)}a_{(s-1)(s-2)}\cdots a_{21})(a_{(t-1)(t-2)}
\cdots a_{(s+1)s})(a_{ts}) $
\end{enumerate}

\item [{\rm (II)}]  Let $A = a_{n_m
n_{m-1}}a_{n_{m-1}n_{m-2}}\cdots  a_{n_2n_1}$
where \ \ $n\geq n_m>n_{m-1}>\cdots >n_1\geq 1$. Then $A$ is
positively  equivalent to a word which begins or ends with $a_{n_t n_s}$,
for any
choice of $n_t, n_s$ with $n\geq n_t > n_s \geq 1$.
\item [{\rm (III)}]  $a_{ts} \delta \doteq \delta a_{(t+1)(s+1)}$, where
subscripts are defined
mod $n$.
\end{enumerate}
\end{lemm}

\proof

(I) With Remark \ref{remark:bandgenerator relations} in mind, choose
any pair of indices $t,s$ with $n\geq t>s\geq 1$. We need to show that
$\delta$ can be represented by a word that begins with $a_{ts}$.
Focus first on the elementary braid $a_{(s+1)s}$ in the expression
for $\delta$ which is given in (\ref{eqn:delta}), and apply the
first of the pair of relations in (\ref{eqn:bandrelation2})
repeatedly to move $a_{(s+1)s}$ to the left (increasing its first
index as you do so) until its name changes to $a_{ts}$. Then apply
the second relation in the pair to move it (without changing its
name) to the extreme left end, vis:
\begin{eqnarray*} \delta
&\doteq& a_{n(n-1)}a_{(n-1)(n-2)}\cdots
a_{(s+2)(s+1)}a_{(s+1)s}a_{s(s-1)}\cdots a_{21}\\ &\doteq&
a_{n(n-1)}a_{(n-1)(n-2)}\cdots a_{(t+1)t}a_{ts}a_{t(t-1)}\cdots
a_{(s+2)(s+1)}a_{s(s-1)}\cdots a_{21}\\ &\doteq&
a_{ts}a_{n(n-1)}\cdots a_{(t+2)(t+1)}a_{(t+1)s}a_{t(t-1)}\cdots
a_{(s+2)(s+1)}a_{s(s-1)}\cdots a_{21}\\
&\doteq&(a_{ts})(a_{n(n-1)}\cdots
a_{(t+2)(t+1)}a_{(t+1)s}a_{s(s-1)}\cdots a_{21})(a_{t(t-1)}\cdots
a_{(s+2)(s+1)})
\end{eqnarray*} We leave it to
the reader to show that the proof works equally well when we move letters
to the
right instead of to the left.

(II) The proof of (II)
is a direct analogy of the proof of (I).

(III) To establish (III), we use (I):
\begin{eqnarray*} &&\delta a_{(t+1)(s+1)}\\
&\doteq& a_{ts}a_{n(n-1)}\cdots
a_{(t+2)(t+1)}a_{(t+1)s}a_{t(t-1)}\cdots a_{(s+2)(s+1)}a_{s(s-1)}\cdots
a_{21}a_{(t+1)(s+1)}\\
&\doteq& a_{ts}a_{n(n-1)}\cdots
a_{(t+2)(t+1)}a_{(t+1)s}a_{(t+1)t}a_{t(t-1)}\cdots
a_{(s+2)(s+1)}a_{s(s-1)}\cdots
a_{21})\\
&\doteq& a_{ts}a_{n(n-1)}\cdots a_{(t+2)(t+1)}a_{(t+1)t}a_{t(t-1)}\cdots
a_{(s+2)(s+1)}a_{(s+1)s}a_{s(s-1)}\cdots a_{21})\\
&\doteq& a_{ts}\delta. ||
\end{eqnarray*}

We move on to the main business of this section, the proof that the
semigroup
$B_n^+$ embeds in $B_n$.  We will use Lemma \ref{lemm:delta} in the
following
way: the inner automorphism defined by conjugation by $\delta$ determines
an
index-shifting automorphism $\tau$ of $B_n$ and $B_n^+$ which is a
useful tool to eliminate repetitious arguments. We define:
$$\tau(a_{ts}) = a_{(t+1)(s+1)} \ \ \ {\rm and}
\ \ \ \tau^{-1}(a_{ts}) = a_{(t-1)(s-1)}.$$

Following the ideas which were first used by Garside \cite{Garside}, the
key step
is to  establish
that there are right and left cancellation laws in $B_n^+$. We remark that
 even though
Garside proved this for Artin's presentation, it does not follow that
it's still true
when one uses the band generator presentation. Indeed, counterexamples were
discovered
by Xu \cite{Xu} and given in \cite{KKL}.

If $X\doteq
Y$ is obtained by a sequence of $t$ single  applications of the defining
relations in (\ref{eqn:bandrelation1}) and (\ref{eqn:bandrelation2}):
$$X\equiv W_0\to W_1\to\cdots\to W_t\equiv Y, $$
then the transformation which takes $X$ to $Y$ will be said to be of
{\em chain-length} $t$.

\begin{thm}[Left ``cancellation'']\label{thm:L}
Let $a_{ts}X\doteq a_{rq}Y$ for some positive words $X,Y$.
Then $X$ and $Y$ are related as follows:
\begin{enumerate}
\item[{\rm (I)}]  If there are only two distinct indices, i.e. $t=r$ and
$s=q$, then
$X\doteq Y$. \item[{\rm (II)}] If there are three distinct indices:
\begin{enumerate}
        \item[{\rm (i)}] If $t=r$ and $q<s$, then $X\doteq a_{sq}Z$ and
        $Y\doteq a_{ts}Z$ for some $Z\in\pos$,
        \item[{\rm (ii)}] If $t=r$ and $s<q$, then $X\doteq a_{tq}Z$ and
        $Y\doteq a_{qs}Z$ for some $Z\in\pos$,
        \item[{\rm (iii)}] If $t=q$, then $X\doteq a_{rs}Z$ and
        $Y\doteq a_{ts}Z$ for some $Z\in\pos$,
        \item[{\rm (iv)}] If $s=r$, then $X\doteq a_{sq}Z$ and
        $Y\doteq a_{tq}Z$ for some $Z\in\pos$,
        \item[{\rm (v)}] If $s=q$ and $r<t$, then $X\doteq a_{tr}Z$ and
        $Y\doteq a_{ts}Z$ for some $Z\in\pos$,
        \item[{\rm (vi)}] If $s=q$ and $t<r$, then $X\doteq a_{rs}Z$ and
        $Y\doteq a_{rt}Z$ for some $Z\in\pos$,
        \end{enumerate}
\item[{\rm (III)}] If the four indices are distinct and if
$(t-r)(t-q)(s-r)(s-q)>0$,
then $X\doteq a_{rq}Z$ and
        $Y\doteq a_{ts}Z$ for some $Z\in\pos$.
\item[{\rm (IV)}] If the four indices are distinct, then:
\begin{enumerate}
        \item[{\rm (i)}] If $q<s<r<t$, then $X\doteq a_{tr}a_{sq}Z$
        and $Y\doteq a_{tq}a_{rs}Z$ for some $Z\in\pos$,
        \item[{\rm (ii)}] If $s<q<t<r$, then $X\doteq a_{tq}a_{rs}Z$
        and $Y\doteq a_{rt}a_{qs}Z$ for some $Z\in\pos$,
        \end{enumerate}
\end{enumerate}
\end{thm}

\proof The proof of the theorem for positive words $X$, $Y$ of word length
$j$
that are positively equivalent via a transformation of chain-length $k$
will be referred to as $T(j,k)$. The proof will be proceeded by an
induction
on $(j,k)$ ordered lexicographically. This induction makes sense because
$T(j,1)$ holds for any $j$. Assume that $T(j,k)$ holds for all
pairs $(j,k)<(l,m)$, that is,
\begin{itemize}
\item[($*$)] $T(j,k)$ is true for $0\le j\le l-1$ and any $k$.
\item[($**$)] $T(l,k)$ is true for $k\le m-1$.
\end{itemize}
Now suppose that $X$ and $Y$ are positive words of length $l$
and $a_{ts}X\doteq a_{rq}Y$ via a transformation of chain-length $m\geq
2$.
Let $a_{\beta\alpha}W$ be the first intermediate word in the sequence of
transformation from $a_{ts}X$ to $a_{rq}Y$.
We can assume that $a_{\beta\alpha}\ne a_{ts}$ and $a_{\beta\alpha}\ne
a_{rs}$,
otherwise we apply the induction hypotheses ($**$) to complete the proof.
Furthermore, since $a_{\beta\alpha}W$ must be obtained from $a_{ts}X$ by
a single application of a defining relation, we see by using ($**$) that
$X\doteq aU$ and $W\doteq bU$ for some distinct generators $a$, $b$ and
a positive word $U$.

For case (I), we see again by using ($**$) that $W\doteq bV$ and
$Y\doteq aV$ for a positive word $V$. Then $W\doteq bU\doteq bV$
implies $U\doteq V$ by ($*$). Thus $X\doteq aU\doteq aV\doteq Y$.

It remains to prove cases II, III, IV (i) and IV(ii). We
fix notation as follows: Using ($**$), $W\doteq BV$ and $Y\doteq AV$ for a
positive
word $V$ and two distinct positive words $A$, $B$ of word length 1 or 2
depending on
$a_{\beta\alpha}$ and $a_{rq}$. When the word length of $A$ and $B$ is 1,
we apply ($*$) to $W\doteq bU\doteq BV$. If $b=B$, $U\doteq V$ and so
$X\doteq aU$ and $Y\doteq AU$ are the required form.
If $b\ne B$, we obtain $U\doteq CQ$ and
$V\doteq DQ$ for a positive word $Q$ and two distinct word $C$ and $D$
of word length 1 or 2 depending on $b$ and $B$. Then we apply some
defining relations to $X\doteq aCQ$ and $Y\doteq ADQ$ to achieve the
desired
forms.

When the word lengths of $A$ and $B$ are 2, we rewrite
$B=b'b''$ in terms of generators. Notice that $b'b''=b''b'$.
Apply ($*$) to either $W\doteq bU\doteq b'(b''V)$ or $W\doteq bU\doteq
b''(b'V)$
to obtain $U\doteq CQ$ and either $b''V\doteq DQ$ or $b'V\doteq DQ$
for a positive word $Q$ and two distinct generators $C$ and $D$. In the
tables below, we will use the symbols so that we have $b''V\doteq DQ$.
Apply ($*$) to $b''V\doteq DQ$. If $b''=D$, we obtain $V\doteq Q$ and then
we apply defining relations to $X\doteq aCV$ and $Y\doteq AV$ to get
the required forms. If $b''\ne D$, we obtain $V\doteq EP$ and $Q\doteq FP$
for
a positive word $P$ and distinct words $E$, $F$ of word length 1 or 2
depending
on $b''$ and $D$. Then we apply some
defining relations to $X\doteq aCFP$ and $Y\doteq AEP$ to achieve the
desired
forms.

The four tables below treat cases II, III, IV(i) and IV(ii). The first
column covers
the possible relative positions of $q,r,s,t,\alpha,\beta$. The second
column contains
one of 4 possible forms $aU,aCQ,aCV,aCFP$ of the word $X$ as explained
 above and
similarly the third column contains one of 4 possible forms
$AU,ADQ,AV,AEP$ of $Y$.
Finally the fourth and fifth columns contain the values of $b$ and $B$,
respectively.

In case (II), it is enough to consider the subcases (i) and (ii) because
the other subcases can be obtained from (i) or (ii) by applying the
automorphism
$\tau$.
But we may also assume that $q<s$, otherwise we switch the roles of $X$
and $Y$.

In case (III), there are actually 4 possible positions of $q,r,s,t$ but
they can be obtained from one position by applying $\tau$. Thus we only
consider case $q<r<s<t$. The table shows all possible cases
required in the induction.

In case (IV), it is again enough to consider the case $q<s<r<t$ and
the table covers all possible inductive steps.

\begin{center}
\small
\let\a\alpha
\let\b\beta
\def\Q#1.#2.#3.#4.#5.{$a_{#1}a_{#2}Q\doteq a_{sq}a_{#3}Q$ &
 $a_{#4}a_{#5}Q\doteq a_{ts}a_{#3}Q$}
\def\V#1.#2.#3.#4.#5.{$a_{#1}a_{#2}V\doteq a_{sq}a_{#3}V$ &
 $a_{#4}a_{#5}V\doteq a_{ts}a_{#3}V$}
\def\P#1.#2.#3.#4.#5.#6.#7.#8.{$\begin{array}{r}
a_{#1}a_{#2}a_{#3}P\\
\doteq a_{sq}a_{#4}a_{#5}P\end{array}$ &
$\begin{array}{r}a_{#6}\b_{#7}a_{#8}P\\
\doteq a_{ts}a_{#4}a_{#5}P\end{array}$}
\def\U{$a_{sq}U$ & $a_{ts}U$}
\def\l#1.#2.{&$a_{#1}$&$a_{#2}$\\ \hline}
\def\L#1.#2.#3.{&$a_{#1}$&$a_{#2}a_{#3}$\\ \hline}
\begin{tabular}{|c|c|c|c|c|}\hline
 & $X$ & $Y$ & $b$ & $B$ \\ \hline
$q<s<t=r<\a<\b$&    \Q\b\a.sq.\b\a.\b\a.ts. \l ts. rq.
$q<s<t=r=\a<\b$&    \Q\b s.sq.\b q.\b q.ts. \l ts. tq.
$q<s=\a<t=r<\b$&    \V\b s.sq.\b q.ts.\b q. \L \b t.\b t.sq.
$q<\a<s<t=r<\b$&    \P\b\a.\b s.\a q.s\a.\b q. t\a.\b q.ts. \L ts.\b t.\a
q.
$q=\a<s<t=r<\b$&    \Q\b\a.\b s.\b\a.\b\a.ts.   \l ts.\b t.
%

$q<s<\a<t=r=\b$&    \Q t\a.sq.t\a.t\a.\a s. \l \a s.\a q.
$q<\a<s<t=r=\b$&    \Q s\a.\a q.s\a.t\a.ts. \l ts.\a q.
$\a<q<s<t=r=\b$&    \Q s\a.sq.q\a.q\a.ts.   \l ts.tq.
$q<s<\a<\b<t=r$&    \Q \b\a.sq.\b\a.\b\a.ts.    \l ts.tq.
$q<s=\a<\b<t=r$&    \Q t\b.sq.t\b.\b s.ts.  \l ts.tq.
$q<\a<s=\b<t=r$&    \Q s\a.\a q.s\a.s\a.t\a.    \l t\a.tq.
$q=\a<s=\b<t=r$&    \U  \l tq.tq.
$\a<q<s=\b<t=r$&    \V s\a.sq.q\a.ts.q\a.   \L t\a.t\a.sq.
$q<\a<\b<s<t=r$&    \Q \b\a.sq.\b\a.\b\a.ts.    \l ts. tq.
$q=\a<\b<s<t=r$&    \Q \b q.sq.s\b.t\b.ts.  \l ts. tq.
$\a<q<\b<s<t=r$&    \P \b\a.\b q.s\a.s\b.q\a.t\b.q\a.ts.    \L ts.\b
q.t\a.
$\a<q=\b<s<t=r$&    \Q q\a.s\a.q\a.q\a.ts.  \l ts. t\a.
\end{tabular}
\end{center}

\newpage

\begin{center}
\small
\let\a\alpha
\let\b\beta
\def\Q#1.#2.#3.#4.#5.{$a_{#1}a_{#2}Q\doteq a_{rq}a_{#3}Q$ &
 $a_{#4}a_{#5}Q\doteq a_{ts}a_{#3}Q$}
\def\P#1.#2.#3.#4.#5.#6.#7.#8.{$\begin{array}{r}
a_{#1}a_{#2}a_{#3}P\\
\doteq a_{rq}a_{#4}a_{#5}P\end{array}$ &
$\begin{array}{r}a_{#6}a_{#7}a_{#8}P\\
\doteq a_{ts}a_{#4}a_{#5}P\end{array}$}
\def\l#1.#2.{&$a_{#1}$&$a_{#2}$\\ \hline}
\def\L#1.#2.#3.{&$a_{#1}$&$a_{#2}a_{#3}$\\ \hline}
\begin{tabular}{|c|c|c|c|c|}\hline
 & $X$ & $Y$ & $b$ & $B$ \\ \hline
$q<r<s<t<\a<\b$&    \Q \b\a.rq.\b\a.\b\a.ts.    \l ts. rq.
$q<r<s<t=\a<\b$&    \Q \b s.rq.\b s.\b t.ts.    \l ts.  rq.
$q<r<s=\a<t<\b$&    \Q \b s.rq.\b s.\b s.\b t.  \l \b t.    rq.
$q<r<\a<s<t<\b$&    \Q \b\a.rq.\b\a.\b\a.ts.    \l ts.rq.
$q<r=\a<s<t<\b$&    \Q \b r.rq.\b q.\b q.ts.    \l ts.rq.
$q<\a<r<s<t<\b$&    \P \b\a.\b r.\a q.\b q.r\a.\b q.r\a.ts. \L ts.\b r.\a
q.
$q=\a<r<s<t<\b$&    \Q \b q.\b r.\b q.\b q.ts.  \l ts.\b r.
$q<r<s<\a<t=\b$&    \Q t\a.rq.t\a. t\a.\a s.    \l \a s. rq.
$q<r<\a<s<t=\b$&    \Q s\a.rq.s\a.t\a.ts.   \l ts. rq.
$q<r=\a<s<t=\b$&    \Q sr.rq.sq.tq.ts.  \l ts.rq.
$q<\a<r<s<t=\b$&    \P s\a.sr.\a q.sq.r\a.tq.r\a.ts.    \L ts.tr.\a q.
$q=\a<r<s<t=\b$&    \Q sq.sr.sq. tq.ts. \l ts.tr.
$\a<q<r<s<t=\b$&    \Q s\a.rq.s\a.t\a.ts.   \l ts.rq.
%
$q<r<s<\a<\b<t$&    \Q \b\a.rq.\b\a.    \b\a.ts.    \l ts.rq.
$q<r<s=\a<\b<t$&    \Q t\b.rq.t\b.\b s.ts.  \l ts. rq.
$q<r<\a<s=\b<t$&    \Q s\a.rq.s\a.s\a.t\a.  \l t\a. rq.
$q<r=\a<s=\b<t$&    \Q sr.rq.sq. sq.tq. \l tr. rq.
$q<\a<r<s=\b<t$&    \P s\a.sr.\a q.sq.r\a. sq.r\a.tq.   \L t\a.sr.\a q.
$q=\a<r<s=\b<t$&    \Q sq.sr.sq.sq.tq.  \l tq.sr.
$\a<q<r<s=\b<t$&    \Q s\a.rq.s\a. s\a.t\a. \l t\a.rq.
$q<r<\a<\b<s<t$&    \Q \b\a..rq.\b\a.\b\a.ts.    \l ts.rq.
$q<r=\a<\b<s<t$&    \Q \b r.rq.\b q.\b q.ts.    \l ts.rq.
$q<\a<r<\b<s<t$&    \P \b\a.\b r.\a q.\b q.r\a.\b q.r\a.ts. \L ts.\b r.\a
q.
$q=\a<r<\b<s<t$&    \Q \b q.\b r.\b q.  \b q.ts.    \l ts.\b r.
$\a<q<r<\b<s<t$&    \Q \b\a.rq.\b\a.\b\a.ts.    \l ts.rq.
$q<\a<r=\b<s<t$&    \Q r\a.\a q.r\a.r\a.ts. \l ts.\a q.
$\a<q<r=\b<s<t$&    \Q r\a.rq.q\a.q\a.ts.   \l ts.rq.
$q<\a<\b<r<s<t$&    \Q \b\a.rq.\b\a.\b\a.ts.    \l ts.rq.
$q=\a<\b<r<s<t$&    \Q \b q.rq.r\b. r\b.ts. \l ts.  rq.
$\a<q<\b<r<s<t$&    \P \b\a.r\a.\b q.r\b.q\a.r\b.q\a.ts.    \L ts.r\a.\b
q.
$\a<q=\b<r<s<t$&    \Q q\a.r\a.q\a.q\a.ts.  \l ts.r\a.
\end{tabular}
\end{center}

\newpage

\begin{center}
\small
\let\a\alpha
\let\b\beta
\arraycolsep 0pt
\def\Q#1.#2.#3.#4.#5.#6.#7.{$\begin{array}{r}a_{#1}(a_{#2}a_{#3})Q\\
\doteq (a_{tr}a_{sq})a_{#4}Q\end{array}$ &
 $\begin{array}{r}a_{#5}(a_{#6}a_{#7})Q\\
 \doteq (a_{rs}a_{tq})a_{#4}Q\end{array}$}
\def\QA{$a_{tr}a_{sq}Q$ & $a_{rs}a_{tq}Q$}
\def\PX#1.#2.#3.#4.#5.#6.{$\begin{array}{r}
 a_{#1}a_{#2}(a_{#3}a_{#4})P\\
 \doteq (a_{tr}a_{sq})a_{#5}a_{#6}P\end{array}$ &}
\def\PY#1.#2.#3.#4.#5.#6.{$\begin{array}{r}a_{#1}a_{#2}a_{#3}a_{#4}P\\
\doteq (a_{rs}a_{il})a_{#5}a_{#6}P\end{array}$}
\def\PA#1.#2.#3.#4.#5.#6.#7.{$\begin{array}{r}a_{#1}a_{#2}a_{#3}P\\
 \doteq (a_{tr}a_{sq})a_{#4}P\end{array}$ &
 $\begin{array}{r}a_{#5}a_{#6}a_{#7}P\\
 \doteq (a_{rs}a_{tq})a_{#4}P\end{array}$}
\def\l#1.#2.{&$a_{#1}$&$a_{#2}$\\ \hline}
\def\L#1.#2.#3.{&$a_{#1}$&$a_{#2}a_{#3}$\\ \hline}
\begin{tabular}{|c|c|c|c|c|} \hline
 & $X$ & $Y$ & $b$ & $B$ \\ \hline
$q<s<r<t<\a<\b$&    \Q\b\a.tr.sq.\b\a.\b\a.tq.rs.   \l ts.rq.
$q<s<r<t=\a<\b$&    \Q\b s.tr.sq.\b q.\b t.tq.rs.   \l ts.rq.
$q<s=\a<r<t<\b$&    \PA\b s.tr.\a q.\b q.\b q.rs.\b t.  \L \b t.\b r.\a q.
$q<\a<s<r<t<\b$&    \PX\b\a.\b s.tr.\a q.\b q.s\a. \PY\b q.r\a.\b t.rs.\b
q.s\a.    \L ts.\b r.\a q.
$q=\a<s<r<t<\b$&    \Q\b q.\b s.tr.\b q.\b q.\b t.rs.   \l ts.\b r.
$q<s<r<\a<t=\b$&    \Q t\a.\a r.sq.t\a.t\a.\a q.rs. \l \a s. rq.
$q<s<r=\a<t=\b$&    \QA \l rs. rq.
$q<s<\a<r<t=\b$&    \PA t\a.tr.sq.r\a.tq.r\a.\a s.  \L \a s.tr. \a q.
$q<\a<s<r<t=\b$&    \PA s\a.tr.\a q.s\a.tq.r\a.rs.  \L ts.tr.\a q.
$q=\a<s<r<t=\b$&    \QA \l ts.tr.
$\a<q<s<r<t=\b$&    \Q s\a.tr.sq.q\a.t\a.\b q.rs.   \l ts.rq.
$q<s<r<\a<\b<t$&    \Q\b\a.tr.sq.\b\a.\b\a.tq.rs.   \l ts.rq.
$q<s<r=\a<\b<t$&    \Q\b r.tr.sq.t\b.\b q.tq.rs.    \l ts.rq.
$q<s<\a<r<\b<t$&    \PX\b\a.\b r.t\a.sq.t\b.r\a. \PY\b q.r\a.tq.\a
s.t\b.r\a.   \L ts.\b r.\a q.
$q<s=\a<r<\b<t$&    \PA t\b.\b r.sq.t\b.\b q.rs.tq. \L ts.\b r.\a q.
$q<s<\a<r=\b<t$&    \Q r\a.t\a.sq.r\a.r\a.tq.\a s.  \l ts. rq.
$q<s=\a<r=\b<t$&    \QA     \l ts.sq.
$q<s<\a<\b<r<t$&    \Q\b\a.tr.sq.\b\a.\b\a.tq.rs.   \l ts.rq.
$q<s=\a<\b<r<t$&    \Q t\b.tr.sq.r\b.\b s.tq.rs.    \l ts.  rq.
$q<\a<s=\b<r<t$&    \Q s\a.tr.\a q.s\a.s\a.tq.r\a.  \l t\a.rq.
$q=\a<s=\b<r<t$&    \QA \l tq. rq.
$\a<q<s=\b<r<t$&    \PA s\a.tr.sq.q\a.rs.q\a.t\a.   \L t\a. rq.sq.
$q<\a<\b<s<r<t$&    \Q\b\a.tr.sq.\b\a.\b\a.rs.tq.   \l ts.rq.
\end{tabular}\\
\begin{tabular}{|c|c|c|c|c|}\hline
 & $X$ & $Y$ & $b$ & $B$ \\ \hline
$q=\a<\b<s<r<t$&    \Q\b q.tr.sq.s\b.r\b.tq.rs. \l ts. rq.
$\a<q<\b<s<r<t$&    \PX\b\a.\b q.tr.s\a.s\b.q\a. \PY
r\b.q\a.rs.t\a.s\b.q\a.    \L ts.\b q.r\a.
$\a<q=\b<s<r<t$&    \Q q\a.tr.s\a.q\a.q\a.t\a.rs.   \l ts. r\a.
\end{tabular}
\end{center}
\noindent This completes the proof that `Left Cancellation' is possible in
the
monoid of positive words. \endpf

Similarly we can prove the following theorem.

\begin{thm}[Right ``cancellation'']\label{thm:R}
Let $Xa_{ts}\doteq Ya_{rq}$ for some positive words $X,Y$.
Then $X$ and $Y$ are related as follows:
\begin{enumerate}
\item[{\rm (I)}] If $t=r$ and $s=q$, then $X\doteq Y$,
\item[{\rm (II)}]\begin{enumerate}
        \item[{\rm (i)}] If $t=r$ and $q<s$, then $X\doteq Za_{tq}$ and
        $Y\doteq Za_{sq}$ for some $Z\in\pos$,
        \item[{\rm(ii)}] If $t=r$ and $s<q$, then $X\doteq Za_{qs}$ and
        $Y\doteq Za_{ts}$ for some $Z\in\pos$,
        \item[{\rm(iii)}] If $t=q$, then $X\doteq Za_{rt}$ and
        $Y\doteq Za_{ts}$ for some $Z\in\pos$,
        \item[{\rm(iv)}] If $s=r$, then $X\doteq Za_{tq}$ and
        $Y\doteq Za_{ts}$ for some $Z\in\pos$,
        \item[{\rm(v)}] If $s=q$ and $r<t$, then $X\doteq Za_{rs}$ and
        $Y\doteq Za_{tr}$ for some $Z\in\pos$,
        \item[{\rm(vi)}] If $s=q$ and $t<r$, then $X\doteq Za_{rt}$ and
        $Y\doteq Za_{ts}$ for some $Z\in\pos$,
        \end{enumerate}
\item[{\rm(III)}] If $(t-r)(t-q)(s-r)(s-q)>0$, then $X\doteq Za_{rq}$ and
        $Y\doteq Za_{ts}$ for some $Z\in\pos$.
\item[{\rm(IV)}]\begin{enumerate}
        \item[{\rm (i)}] If $q<s<r<t$, then $X\doteq Za_{tq}a_{rs}$
        and $Y\doteq Za_{tr}a_{sq}$ for some $Z\in\pos$,
        \item[{\rm (ii)}] If $s<q<t<r$, then $X\doteq Za_{rt}a_{qs}$
        and $Y\doteq Za_{tq}a_{rs}$ for some $Z\in\pos$,
        \end{enumerate}
\end{enumerate}
\end{thm}

The properties of $\delta$ which were worked out in Lemma \ref{lemm:delta}
ensure that
$\delta$ can take the role of the half twist $\Delta$ of the Garside's
argument
in~\cite{Garside} to show:

\begin{thm}[Right reversibility]\label{thm:reversibility}
If $X,Y$ are positive words,
then there exist positive words $U,V$ such that $UX\doteq VY$.
\end{thm}

Using left and right ``cancellation'' and right
reversibility, we obtain (as did Garside) the following embedding
theorem~\cite{Clifford:Preston}.

\begin{thm}[Embedding Theorem]\label{thm:embedding}
The natural map from $\pos$ to $B_n$ is injective, that is,
if two positive words are equal in $B_n$,
then they are positively equivalent.
\end{thm}

\begin{remark}
{\rm Any time that the defining relations
in a group presentation are expressed as relations between positive words in
the generators one may consider the semigroup of positive words and ask
whether that semigroup embeds in the corresponding group.  Adjan
\cite{Adjan} and also Remmers \cite{Remmers} studied this
situation and showed that a semigroup is embeddable if it is
`cycle-free', in their terminology. Roughly speaking, this means that
the presentation has relatively few relations, so that a positive word
can only be written in a small number of ways. But the fundamental
words $\Delta$ and $\delta$ can be written in many many ways, and it
therefore follows that large subwords of these words can too, so our
presentations are almost the opposite to those considered by Adjan and
Remmers.

    According to Sergiescu \cite{Sergiescu}, any connected
planar graph
with $n$ vertices gives rise to a positive presentation of
$B_n$ in
which each edge gives a generator which is a conjugate of
one of Artin's
elementary braids and relations are derived at each vertex
and at each
face. In fact one can generalize his construction as follows.
Consider the
elements in $B_n$ defined by:
$$b_{ts}=(\sigma_{t-1}^{-1}\sigma_{t-2}^{-1}\cdots
\sigma_{s+1}^{-1})\sigma_s
(\sigma_{s+1}\cdots\sigma_{t-2}\sigma_{t-1}).$$
The braid $b_{ts}$ is geometrically a positive half-twisted band
connecting the $t^{th}$
and the $s^{th}$ strands, and passing  behind all intermediate strands.
Since $b_{t(t-1)}=\sigma_i=a_{t(t-1)}$,
the set $X=\{a_{ts}, b_{ts}\vert 1\le s< t\le n\}$ contains
$(n-1)^2$
elements. Then $X$ may be described by a graph on a plane where
$n$ vertices
are arranged, in order, on a line. An element $a_{ts}, b_{ts} \in X$
belongs to an edge connecting the $t^{th}$ and the
$s^{th}$ vertices
on one side or the other, depending on whether it is an $a_{ts}$ or
$b_{ts}$. In this way one obtains a planar graph
in which two edges have at most one interior intersection point. It
is not hard
to show that a subset $Y\subset X$ is a generating set of $B_n$
if and only
if the generators in $Y$ form a connected subgraph. Consider all
presentations that have $Y \subset X$ as a set of
generators and have a finite set of equations between positive
words in $Y$
as a set of relators. All Sergiescu's planar graphs are of this type.
Artin's presentation corresponds to the linear graph with $n-1$ edges
and our presentation corresponds to the complete graph on $n$ vertices.
One can prove \cite{Ko} that the embedding theorem fails to
hold in all but two presentations of this type. Those two are
Artin's presentation and ours. }
\end{remark}

\newpage

\section{The word problem}

In this section we present our solution to the word problem in $B_n$,
using the presentation of Proposition \ref{thm:presentation}.  Our
approach builds on the ideas of Garside \cite{Garside},  Thurston
\cite{Epstein} and Elrifai and Morton in
\cite{Elrifai:Morton}.  In the next section we will translate the results of
this section into an algorithm, and compute its complexity.

We begin with a very simple consequence of Lemma \ref{lemm:delta}.

\begin{lemm}
\label{positive and negative}
Every element ${\cal W}\in B_n$ can be represented by a word of the form
$\delta^pQ$ where $p$ is an integer and $Q$ is a word in the
generators $a_{t,s}$ of $B_n^+\subset B_n$.
\end{lemm}

\proof Choose any word which represents ${\cal W}$. Using (I)
of Lemma \ref{lemm:delta} replace every generator which occurs with a
negative exponent by $\delta^{-1}M$, where $M$ is positive. Then use
(III) of Lemma \ref{lemm:delta} to collect the factors $\delta^{-1}$ at
the left. \endpf

The {\it word length} of a (freely reduced) word $W$ in our presentation
of
$B_n$ is denoted by $|W|$. The identity
word will be denoted by
$e$,
$|e|=0$. For words $V,W$, we write $V\le W$ (or
$W\geq V$) if
$W=P_1VP_2$ for some $P_1,P_2\in\pos$.
Then $W\in\pos$ if and only if $e\le W$. Also
$V\le W$ if and only if $W^{-1}\le V^{-1}$.

Recall that $\tau$ is the inner automorphism of $B_n$ which is defined by
$\tau(W)=\delta^{-1}W\delta$. By Lemma \ref{lemm:delta} the
action of $\tau$ on the generators is given by $\tau(a_{ts}) =
a_{(t+1)(s+1)}$.

\begin{prop}
\label{prop;partial order}
The relation `$\le$' has the following properties:
\begin{enumerate}
\item [{\rm (I)}] `$\le$' is a partial order on $B_n$.
\item [{\rm (II)}] If $W\le\delta^u$, then $\delta^u=PW=W\tau^u(P)$ for some
    $P\in\pos$
\item [{\rm (III)}] If $\delta^u\le W$, then $W=P\delta^u=\delta^u\tau^u(P)$
    for some $P\in\pos$
\item [{\rm (IV)}] If $\delta^{v_1}\le V\le\delta^{v_2}$ and
    $\delta^{w_1}\le W\le\delta^{w_2}$, then
    $\delta^{v_1+w_1}\le VW\le\delta^{v_2+w_2}$.
\item  [{\rm (V)}] For any $W$ there exist integers $u,v$
such that  $\delta^u\le W\le\delta^v$.
\end{enumerate}
\end{prop}

\proof See \cite{Elrifai:Morton} or \cite{KKL}. The proofs given
there carry over without any real changes to the new situation. \endpf

\

The set $\{W\mid\delta^u\le W\le\delta^v\}$ is denoted by $[u,v]$.
For ${\cal W}\in B_n$, the last assertion of the previous proposition
enables us to define the {\it infimum} and the {\it supremum} of ${\cal W}$
as
$\inf ({\cal W}) =\max\{u\in\ints\mid\delta^u\le W\}$ and
$\sup ({\cal W})=\min\{v\in\ints\mid W\le\delta^v\}$, where $W$ represents
${\cal W}$. The integer
$\ell ({\cal W})=\sup({\cal W})-\inf({\cal W})$ is called the {\em canonical
length} of
${\cal W}$.

A permutation $\pi$ on $\{1,2,\ldots,n\}$ is called a {\em descending
cycle}
if it is represented by a cycle
$(t_j,t_{j-1},\ldots,t_1)$ with $j\geq 2$ and $t_j>t_{j-1}>\ldots>t_1$.
Given a descending cycle
$\pi=(t_j,t_{j-1},\ldots,t_1)$, the symbol $\delta_\pi$ denotes the
positive
braid $a_{t_jt_{j-1}}a_{t_{j-1}t_{j-2}}\cdots a_{t_2t_1}$.
A pair of descending cycles $(t_j,t_{j-1},\ldots,t_1)$,
$(s_i,s_{i-1},\ldots,s_1)$  are said to be {\em parallel} if $t_a$ and
$t_b$
never separate $s_c$ and $s_d$.  That is,
$(t_a-s_c)(t_a-s_d)(t_b-s_c)(t_b-s_d)>0$
for all $a,b,c,d$ with $1\le a<b\le j$ and $1\le c<d\le i$. The cycles in
a product of
parallel, descending cycles are disjoint and non-interlacing.
Therefore they commute with one-another.
For pairwise parallel,
descending cycles $\pi_1$, $\pi_2$,$\ldots$, $\pi_k$,
the factors in the product
$\delta_{\pi_1}\delta_{\pi_2}\cdots\delta_{\pi_k}$
are
positive braids which commute with one-another and therefore there is
a well-defined map from the set of all products of parallel descending
cycles to
$B_n$, which splits the homomorphism $\phi:B_n\to\Sigma_n, \ \ \phi(a_{ts})
=(t,s)$.

Our first goal is to prove that braids in $[0,1]$, i.e. braids $A$
with $e \leq A \leq \delta$ are precisely the products
$\delta_{\pi_1}\delta_{\pi_2}\cdots\delta_{\pi_k}$ as above. We will also
prove
that each $\delta_{\pi_i}$ is represented by a unique word in the band
generators, and
so that the product $A$ also has a representation which is unique up
to the order of the factors.

Let $A = BaCbD$ be a decomposition of the positive word $A$ into subwords,
where
$a,b$ are generators.  Let $t,s,r,q$ be integers, with  $n\geq t > s > r >
q
\geq 1$.
We say that the pair of letters $(a,b)$ is {\em an obstructing pair} in the
following cases:

\

case (1): \ $a=a_{tr}, \ b=a_{sq}$

case (2): \ $a=a_{sq}, \ b=a_{tr}$

case (3): \ $a=a_{sr}, \ b=a_{ts}$

case (4): \ $a=a_{ts}, \ b=a_{tr}$

case (5): \ $a=a_{tr}, \ b=a_{sr}$

case (6): \ $a=a_{ts}, \ b=a_{ts}$.

\begin{lemm}\label{lemm:BaCbD}
A necessary condition for  a positive word $A$ to be in $[0,1]$ is that
$A$ has
no  decomposition as $BaCbD$, with $B,a,C,b,D \geq e$ and $(a,b)$ an
obstructing pair.\end{lemm}

\proof
We use a geometric argument.
Given a braid word $W$ in the $a_{t,s}$'s, we associate to $W$ a surface
$F_W$
bounded  by the closure of $W$, as follows: $F_W$ consists of $n$
disks joined by half-twisted bands, with a band for each letter in $W$.
The
half-twisted band for $a_{ts}$ is the negative band connecting
the
$t^{th}$ and the $s^{th}$ disks. Our defining relations in
(\ref{eqn:bandrelation1}) and (\ref{eqn:bandrelation2}) correspond to
isotopies sliding a half-twisted band over an adjacent half-twisted band
or moving a half-twisted band horizontally. (See
Figure 1(b)).
Thus defining relations preserve the topological characteristics of $F_W$.
For
example the surface $F_\delta$ has one connected component and is
contractible.

By the proof of Lemma~\ref{lemm:delta}, we may write $\delta=a_{ts}W$
where
$W=\delta_{\pi'}\delta_{\pi''}$ for parallel descending cycles
$\pi'=(t,t-1,\ldots,s+1)$ and $\pi''=(n,n-1,\ldots,t+1,s,s-1,\ldots,1)$.
Thus
for this $W$ the surface $F_W$ has two connected components,
$F_{\delta_{\pi'}}$
and $F_{\delta_{\pi''}}$.

It is enough to consider the cases $(a,b) = (a_{tr}, a_{sq}),
(a_{sr}, a_{ts}), (a_{ts},a_{ts})$
since all other cases are obtained from these cases by applying the
automorphism
$\tau$, which preserves $\delta$.  Since $A$ is in $[0,1]$ we know that
$\delta = V_1AV_2$ for some $V_1,V_2\geq e$. By Proposition
\ref{prop;partial
order} (II) we see that $AE=\delta$ for some word $E\geq e$. So $BaCbDE
= \delta$, which implies that $aCbDE\tau(B)=\delta$. If $a=a_{tr},
b=a_{sq}$ for
$t>s>r>q$ and $\delta=a_{tr}W$, then $F_W$ has two connected components
and the
$s^{th}$ disk and the $q^{th}$ disk lie on distinct components. But the
$s^{th}$ and
the $q^{th}$ disks lie in the same component in $F_{CbDE\tau(B)}$ since
they are
connected by $b$ and this is a contradiction.

If $a=a_{sr}, b=a_{ts}$ and $\delta=a_{sr}W$, then the $t^{th}$ and
$s^{th}$ disks lie
in distinct components in $F_W$ but they lie in the same component in
$F_{CbDE\tau(B)}$
and this is again a contradiction.

If $a=b$, then $F_{AE}$ contains a non-trivial loop but
$F_\delta$ is contractible and this is a contradiction. \endpf

\begin{thm}\label{thm:subwords}
A braid word $A$ is in $[0,1]$ if and only if
$A=\delta_{\pi_1}\delta_{\pi_2}\cdots\delta_{\pi_k}$
for some parallel, descending cycles $\pi_1$, $\pi_2$,$\ldots$, $\pi_k$ in
$\Sigma_n$.
\end{thm}

\proof  First assume that
$A=\delta_{\pi_1}\delta_{\pi_2}\cdots\delta_{\pi_k}$.
We induct on the number $n$ of braid strands to
prove the necessity. The theorem is true when $n=2$.
Suppose that $\pi_1$, $\pi_2$,\ldots,$\pi_k$ are parallel, descending
cycles in
$\Sigma_n$. In view of the inductive hypothesis, we may assume without
loss of
generality that the index $n$ appears in one of cycles. Since
the factors $\delta_{\pi_1},\dots,\delta_{\pi_j}$ in the product commute
with
one-another, we may assume that for some $1\le i\le k$, the cycle
$\pi_i=(n,t,\ldots,s)$, where all of the indices occurring in
$\pi_1,\ldots,\pi_{i-1}$ are greater than $t$ and all of the indices
occurring in
$\pi_{i+1},\ldots,\pi_k$ are less than $t$. The induction hypothesis
implies that \begin{eqnarray*}
C_1\delta_{\pi_1}\cdots\delta_{\pi_{i-1}}&=&a_{(n-1)(n-2)}\cdots
a_{(t+1)t}\\
\delta_{\pi_i^\prime}\delta_{\pi_{i+1}}\cdots\delta_{\pi_k}C_2
&=&a_{t(t-1)}\cdots a_{21} \end{eqnarray*} where
$\pi_i^\prime=(t,\ldots,s)$ and
$C_1$, $C_2$ are positive words. Thus $$
C_1AC_2=a_{(n-1)(n-2)}\cdots a_{(t+1)t}a_{nt}a_{t(t-1)}\cdots a_{21}=
a_{n(n-1)}a_{(n-1)(n-2)}\cdots a_{21} = \delta. $$
Thus our condition is necessary.

Now assume that $A$ is in $[0,1]$. We prove sufficiency by induction on
the word length of $A$. The theorem is true when $|A|=1$.
Suppose, then, that $|A|>1$.   Let $A=a_{ts}A'$.
By the induction hypothesis $A'=\delta_{\pi_1}\delta_{\pi_2}\cdots
\delta_{\pi_k}$ for some parallel, descending cycles $\pi_1$,
$\pi_2$,$\ldots$, $\pi_k$ in $\Sigma_n$. Since $A$ is in $[0,1]$, we know,
from Lemma \ref{lemm:BaCbD}, that $A$ has no decomposition as $BaCbD$ with
$(a,b)$ an obstructing pair, so in particular there is no $a_{rq}\in
A^\prime$ such that
$(a_{ts},a_{rq})$ is an obstructing pair. Therefore, in particular, by
cases (1) and (2) for obstructing pairs we must have
$(t-r)(t-q)(s-r)(s-q)\geq 0$ for all $a_{rq}\le A'$.
Therefore, if neither $t$ nor $s$ appears among the indices in
any of the $\pi_i$,
then the descending cycle $(t,s)$ is clearly parallel to each $\pi_i$ and
$A=a_{ts}\pi_1\cdots\pi_k$ is in the desired form.

Suppose that $t$ appears in some $\pi_i=(t_1,t_2,\ldots,t_m)$.
Then, by cases (3) and (4) for obstructing
pairs we must have
$t=t_1$ and $s<t_m$.
Suppose that $s$ appears in some $\pi_j=(s_1,s_2,\ldots,s_l)$.
Then case (5) in our list of obstructing pairs tells us that either
$s=s_h$ and
$t<s_{h-1}$
for $1<h\le l$ or $s=s_1$.
Thus we have the following three possibilities:
\begin{enumerate}
\item[(i)] $t$ appears in some $\pi_i=(t_1,t_2,\ldots,t_m)$ and $s$ does
not appear.
Then
$$A=\delta_{\pi_1}\cdots\delta_{\pi_{i-1}}\delta_{\pi'_i}
\delta_{\pi_{i+1}}\cdots
\delta_{\pi_k}$$
is in the desired form, where $\pi'_i=(t_1,t_2,\ldots,t_m,s)$;
\item[(ii)] $s$ appears in some $\pi_j=(s_1,s_2,\ldots,s_l)$ and $t$ does
not appear.
Then
$$A=\delta_{\pi_1}\cdots\delta_{\pi_{j-1}}\delta_{\pi'_j}
\delta_{\pi_{j+1}}\cdots
\delta_{\pi_k}$$ is in the desired form, where
$\pi'_j=(s_1,\ldots,s_{h-1},t,s_h,\ldots,s_l)$ or $(t,s_1,\ldots,$ $s_l)$;
\item[(iii)] $t$ appears in some $\pi_i=(t_1,t_2,\ldots,t_m)$ and
$s$ appears in some $\pi_j=(s_1,s_2,\ldots,s_l)$. Then we may assume $i<j$
and
$$A=\delta_{\pi_1}\cdots\delta_{\pi_{i-1}}\delta_{\pi'_i}
\delta_{\pi_{i+1}}\cdots
\delta_{\pi_{j-1}}\delta_{\pi_{j+1}}\cdots\delta_{\pi_k}$$
is in the desired form, where
$\pi'_i=(t_1,t_2,\ldots,t_m,s_1,s_2,\ldots,s_l)$. $||$
\end{enumerate}

\

\noindent {\bf Definition.} \ From now on we will refer to a braid
which is in
$[0,1]$, and which can therefore be represented by a product of parallel
descending cycles, as a {\em canonical factor}.  For example,
the
14 distinct canonical factors for $n=4$ are:
$$e, a_{21}, a_{32},
a_{31}, a_{43}, a_{42}, a_{41}, a_{32}a_{21},  a_{43}a_{32}, a_{43}a_{31},
a_{43}a_{21},  a_{42}a_{21}, a_{41}a_{32}, a_{43}a_{32}a_{21}.$$
A somewhat simpler notation
describes a descending cycle by its subscript array.
In the example just given the 13 non-trivial canonical factors are:
$$(21),
(32), (31), (43), (42), (41), (321), (432), (431), (421), (4321), (43)(21),
(41)(32).$$  The associated permutation is the cycle associated to
the reverse of the subscript array, with all
indices which are not listed explicitely fixed.

\begin{coro}
\label{cor;Catalan number} For each fixed positive integer $n$ the
number of distinct canonical factors is the $n^{th}$ Catalan number
${\cal C}_n = {(2n)!}/{n!(n+1)!}$.
\end{coro}

\proof
We associate to each product
$\pi=\pi_1\pi_2\cdots\pi_k$ of parallel descending cycles a set of $n$
disjoint arcs in the upper half-plane whose $2n$ endpoints are on the real
axis. Mark the numbers $1,2,\dots, n$ on the real axis. Join $i$ to $\pi(i)$
by an
arc, to obtain $n$ arcs, some of which may be loops. Our arcs have
disjoint
interiors because the cycles in $\pi$ are parallel. By construction there
are
exactly  two arcs meeting at each integer point on the real axis. Now
split
the $i^{th}$ endpoint, $i=1,2,\dots,n$, into two points, $i^\prime,
i^{\prime\prime}$, to obtain $n$ disjoint arcs with
$2n$ endpoints. The pattern so obtained will be called an
[n]-configuration.  To recover the product of disjoint cycles, contract
each
interval
$[i^\prime,i^{\prime\prime}]$ to a single point $i$. In this way we see
that
there is a one-to-one correspondence between canonical factors and
[n]-configurations. But the number of [n]-configurations is the
$n^{th}$ Catalan number (see \cite{KS} for a proof).
\endpf

Note that ${\cal C}_n/{\cal C}_{n-1}=4-\frac6{n+1}\le 4$ and so
${\cal C}_n\le 4^n$. In the Artin presentation of $B_n$, the number of
permutation braids is $n!$ which is much greater than ${\cal C}_n$. This
is one of the reasons why our presentation gives faster algorithm.

It is very easy to recognize canonical factors when they are given as
products of parallel descending cycles. If, however, such a
representative is modified in some way by the defining relations, we
will also need to be able to recognize it.  For computational
purposes the following alternative characterization of
canonical factors will be extremely useful.
It rests on Lemma \ref{lemm:BaCbD}:

\begin{coro}
\label{coro:BaCbD}
A positive word $A$ is a canonical factor if and only if $A$ contains
no obstructing pairs.
\end{coro}

\proof We established necessity in Lemma
\ref{lemm:BaCbD}. We leave it to the reader to check that the proof of
Theorem \ref{thm:subwords} is essentially a proof of sufficiency.
\endpf

The braid $\delta$ can be written in many different ways as a product of
the $a_{ts}$, and by Corollary \ref{coro:BaCbD} each such
product contains no obstructing pair. Any descending cycle
$\delta_\pi$ also has this property. If an element in
$B_n$ is  represented by a word which contains no obstructing
pairs, then it is a canonical factor and so
it can be written as a product of parallel descending cycles.
It follows that there is no obstructing pair in any word representing it.

To get more detailed information about canonical factors
$\delta_{\pi_1}\delta_{\pi_2}\cdots\delta_{\pi_k}$, we begin to
investigate some of their very nice properties.
We proceed as in the
foundational paper of Garside \cite{Garside} and define the {\em
starting set} $S(P)$ and the {\em finishing set} $F(P)$:
\begin{eqnarray*} S(P)&=&\{a\mid P= a P^\prime, P^\prime\geq e, \ \
a\rm{\ \  is \ \  a \ \
 generator}\},\\ F(P)&=&\{a\mid P= P^\prime a, P^\prime\geq e, \ \ a\rm{ \
\ is
\ \ a \ \  generator}\}. \end{eqnarray*}
Note that
$S(\tau(P))=\tau(S(P))$ and $F(\tau(P))=\tau(F(P))$.

Starting sets play a fundamental role in the solutions to the word and
conjugacy problems in \cite{Elrifai:Morton}. Our canonical form
allows us to determine them by inspection.

\begin{coro}
\label{coro:starting sets}
The starting sets of canonical factors satisfy the following
properties: \begin{enumerate}
\item [{\rm (I)}] If $\pi=(t_m,t_{m-1},\ldots,t_1)$ is a descending
cycle,  then the starting set (and also the finishing set)
of $\delta_{\pi}$ is
$\{a_{t_jt_i}; \ m\geq j > i\ \geq 1 \}$.
\item [{\rm (II)}] If $\pi_1,\ldots,\pi_k$ are parallel
descending cycles, then   $S(\delta_{\pi_1}\cdots\delta_{\pi_k})
=S(\delta_{\pi_1})\cup\cdots\cup S(\delta_{\pi_k}).$
\item [{\rm (III)}] If $A$ is a canonical factor,
then $S(A)=F(A)$.
\item[{\rm (IV)}] If $A$ and $B$ are canonical factors, and if $S(A)=S(B)$
then $A=B$.
\item[{\rm (V)}] Let $P$ be a given positive word. Then there exists a
canonical factor $A$ such that $S(P)=S(A)$.
\item[{\rm (VI)}] If
$S(A)\subset S(P)$ for some canonical factor $A$, then $P=AP'$ for some
$P'\geq e$.
\item [{\rm (VII)}] For any $P\geq e$, there is a unique
canonical  factor $A$ such that
$P=AP'$ for some $P'\geq e$ and $S(P)=S(A)$.
\end{enumerate}
\end{coro}

\proof To prove (I), observe that the defining relations (7) and (8)
preserve the set of distinct subscripts which occur in a positive word, so
if $a_{qp}$ is in the starting set (resp. finishing set) of $\delta_\pi$
then $q=t_j$ and $p=t_i$ for some $j,i$ with $m\geq j>i\geq 1$. Since
it is proved in part (II) of Lemma~\ref{lemm:delta} that every
$a_{t_jt_i}$ occurs in both the starting set and the finishing set, the
assertion follows.

To prove (II) one need only notice that $\delta_{\pi_r}$ commutes
with $\delta_{\pi_s}$ when the cycles $\pi_r,\pi_s$ are parallel.

Clearly (III) is a consequence of (I) and (II).

As for (IV), by Theorem~\ref{thm:subwords} a canonical word is uniquely
determined by a set of parallel descending cycles. If two distinct
descending cycles $\pi,\mu$ are parallel, then $\delta_\pi,\delta_\mu$
have distinct starting sets, so if $A$ and $B$ are canonical factors,
with $S(A)=S(B)$, the only possibility is that $A=B$.

To prove (V), we induct on the braid index $n$. The claim
is clear for $n=2$. Let $P\in B_n^+$ have starting set $S(P)$. If all
generators of the form $a_{nt}$ for $n-1\geq t\geq 1$ are deleted from
$S(P)$ we obtain a set $S^\prime(P)$ which, by the induction hypothesis,
is the starting set of a braid
$A^\prime = \delta_{\pi_1}\cdots\delta_{\pi_k}$, where
$\pi_1,\ldots,\pi_k$ are parallel, descending cycles in $\Sigma_{n-1}$.
It is now enough to check the following properties of $S(P)$:
\begin{enumerate}
\item[(i)] If $a_{ns},a_{tr}\in S(P)$, with $t>s$, then $a_{ts}\in
S(P)$;
\item[(ii)] If $a_{ns}\in S(P)$ and if $s$ happens to be in one
of the descending cycles
$\pi_i=(t_m,\ldots,t_1)$ associated to $S^\prime(P)$, then $a_{nt_j}\in
S(P)$ for every $j$ with $m\geq j\geq 1$;
\item[(iii)] If $a_{ns}\in S(P)$, where $s$ is not in any
of the descending cycles $\pi_1^\prime,\dots,\pi_r^\prime$ associated to
$S^\prime(P)$, then there is no $a_{tr}\in S(P)$ such that $t>s>r$.
\end{enumerate}

To establish (i), note that since $a_{ns},a_{tr}\in S(P)$, we
have  $P=a_{ns}X = a_{tr}Y$ with $n>t>s>r$.  But then the assertion
follows from Theorem \ref{thm:L}, part (IV), case (ii).

To establish (ii), set $\mu=(n,t_m,\ldots,t_1).$  Then
$\mu$ is a descending cycle for $A$, so by (I) and (II) of this
lemma we conclude  that $a_{nt_j}\in S(P)$ for $m\geq j\geq 1$.

Property (iii) can be verified by observing that if $a_{ns}\in S(P)$, then
$\mu=(n,\dots,s\dots)$ (where $\mu$ could be $(n,s)$ )
must be a descending cycle belonging to a canonical factor
$\delta_\mu$ for $A$. But if so, and if
$a_{tr}\in S(P)$ with
$r<s<t$, then by (i) $a_{ts}$ is also in $S(P)$, so that in fact $\mu =
(n,\dots,t,\dots,s,\dots,r,\dots)$. But then the cycle
survives after deleting $n$,
contradicting the hypothesis that the subscript $s$ does not appear in any
descending cycle associated to $S^\prime(P)$.  Thus we have proved (V).

To prove (VI), induct on the word length of $P$. The assertion is clear if
$|P|=1$. Assume $|P|>1$. We may assume that $a_{nt}$ is in $S(A)$ for some
$1\le t<n$, otherwise we apply the index-shifting automorphism $\tau$.
We make this assumption to reduce
the number of the cases that we have to consider.
By Theorem \ref{thm:subwords}, we may write
$A=\delta_{\pi_1}\delta_{\pi_2}\cdots\delta_{\pi_k}$
for some parallel, descending cycles $\pi_1$, $\pi_2$,$\ldots$, $\pi_k$ in
$\Sigma_n$
and we may assume $\pi_1=(n,t_1,\ldots,t_j)$.
Let $A=a_{nt_1}B$ and $P=a_{nt_1}Q$. We are done by induction if we show
$S(B)\subset S(Q)$. Let $a_{sr}$ be any member of $S(B)$. We have three
possible
cases after considering the properties of the words $A$ and $B$:
\begin{enumerate}
\item[(i)] $s=t_1$ and $r=t_i$ for some $2\le i\le
j$; \item[(ii)] $s<t_1$;
\item[(iii)] $t_1<r<s<n$.
\end{enumerate}
When (i) is the case, then $A = a_{nt_1}a_{t_1t_i}C =
a_{nt_i}a_{t_1t_i}C$ for some
 canonical factor $C$ and so $a_{nt_i}\in S(P)$. Since both
$a_{nt_1}$ and $a_{nt_i}$ are in $S(P)$, Theorem \ref{thm:L}(II)(ii)
implies that $a_{t_1t_i}\in S(Q)$. For the other two cases we can
show, in a similar way, that $a_{sr}$ is in $S(Q)$, using Theorem
\ref{thm:L}.

Assertion (VII) is an immediate consequence of (V) and (VI).
\endpf

Theorem \ref{thm:subwords} has given us an excellent
description of the canonical factors. What remains is to translate
it into a solution to the word problem.  For that purpose we need
to consider {\em products} $A_1A_2\cdots A_k$, where each $A_i$ is
a canonical factor. The argument we shall use is very similar to that
in  \cite{Elrifai:Morton} and \cite{KKL}, even though our $\delta$ and
our canonical factors are very different from their
$\Delta$ and their permutation braids.

A decomposition $Q=AP$,  where $A$ is a canonical factor and
$P \geq e,$ is said to be {\em left-weighted} if
$|A|$ is maximal for all such decompositions. Notice that $AP$ is
{\it not} left-weighted of there exists $p\in S(P)$ such that $Ap$ is a
canonical factor, for if so then $|A|$ is not maximal. We call
$A$  the {\em maximal head} of $Q$ when $Q=AP$ is left-weighted. The
symbol $A\lceil P$ means that $AP$ is left-weighted.
The following corollary gives an easy way to check whether a given
decomposition is left-weighted.

\begin{coro}\label{coro:lwtest}
Let $A,P$ be positive words, with $A$ representing a canonical factor. Then
$A\lceil P$ if and only if for each $b\in S(P)$ there exists $a\le
A$ such that $(a,b)$ is an obstructing pair.
\end{coro}

\proof By the definition of left-weightedness,
$A\lceil P$ if and only if, for each $b\in S(P)$, $Ab$ is not a canonical
factor. By Corollary \ref{coro:BaCbD} $Ab$ is not a
canonical factor if and only if $Ab$ contains an
obstructing pair $(a,q)$. We cannot have both $a\le A$ and $q\le A$ because
by hypothesis $A$ is a canonical factor so that by Corollary
\ref{coro:BaCbD} no word which represents it contains an obstructing pair.
Thus $q=b$.
\endpf

Define the {\em
right complementary set} $R(A)$  and  the {\em left complementary set}
$L(A)$ of a canonical factor $A$ as follows:
$$ R(A) = \{a \ | \ Aa\leq\delta\}$$
$$ L(A) = \{a \ | \ aA\leq\delta\},$$ where $a$ is a generator.
Define the {\em right complement} of a canonical
factor $A$ to be the word $\ov{A}$ such that $A\ov{A} = \delta$.
Since
$W\delta=\delta\tau(W)$, we have $L(\tau(A))=\tau(L(A))$ and
$R(\tau(A))=\tau(R(A))$.

Note that $\ov{(\ov{A})} = \tau(A), \ \
R(A)=S(\ov{A})=F(\ov{A}), \ \ L(\ov{A})=F(A)=S(A)$. Also
$R(\ov{A})=S(\tau(A))=F(\tau(A))$, because
$\tau^{-1}(\ov{A})A=\delta=\ov{A}\tau(A)$ and
$L(A)=F(\tau^{-1}(\ov{A}))=S(\tau^{-1}(\ov{A}))$.

The next proposition shows us equivalent ways to recognize
when a decomposition of a positive word is left-weighted.

\begin{prop}
\label{prop:TFAE}
For any $Q\geq e$, let $Q=AP$ be a decomposition, where
$A$ is a canonical factor and $P\geq e$. Then the following are
equivalent: \begin{enumerate}
\item[{\rm (I)}] $A\lceil P$.
\item[{\rm (II)}] $R(A)\cap S(P)=\emptyset.$
\item[{\rm (III)}] $S(Q)=S(A)$.
\item[{\rm (IV)}] If $WQ\geq\delta$ for some $W\geq e$, then
$WA\ge\delta$.
\item[{\rm (V)}] For any $V\geq e$, $S(VQ)=S(VA)$.
\item[{\rm (VI)}] If $Q=A_1P_1$ is another decomposition with
$A_1$ a canonical factor and $P_1\geq e$, then $A=A_1A^\prime$
for some canonical factor $A^\prime$ (where $A^\prime$ could be $e$).
\end{enumerate}
\end{prop}

\proof See \cite{Elrifai:Morton} and \cite{KKL}.  \endpf

We can now give the promised normal form, which  solves the word
problem for our new presentation for $B_n$:

\begin{thm}
\label{thm:l canonical form}
Any n-braid ${\cal W}$ has a unique representative $W$
left-canonical form:
\begin{displaymath} W=\delta^u A_1A_2\cdots A_k,
\end{displaymath}
where each adjacent pair $A_iA_{i+1}$ is left
weighted and each $A_i$ is a canonical factor.  In
this representation
$\inf ({\cal W})=u$ and $\sup ({\cal W})=u+k$.
\end{thm}

\proof  For any $W$ representing ${\cal W}$
we first write $W=\delta^vP$ for some positive word $P$ and
a possibly negative integer $v$.  For any $P\geq e$, we then iterate
the left-weighted decomposition $P=A_1P_1,\;\;P_1=A_2P_2,\;\;\ldots $
to obtain $W=\delta^u A_1A_2\cdots A_k,$
where $e< A_i<\delta$ and
$R(A_i)\cap S(A_{i+1})=\emptyset$.
This decomposition is unique, by Corollary
\ref{coro:starting sets}, because  $S(A_iA_{i+1}\cdots
A_k)=S(A_i)$ for $1\le i\le k$.
\endpf

The decomposition of Theorem \ref{thm:l canonical form}
will be called the {\em left-canonical form} of ${\cal W}$.  For future
use, we note one of its symmetries:

\begin{prop}\label{prop:invform}
\begin{enumerate}
\item [{\rm (I)}] Let $A,B$ be canonical factors. Then $A\lceil B$ if and
only if
$\ov{B}\lceil\tau(\ov{A})$.
\item [{\rm (II)}] The left-canonical forms of ${\cal W}$ and
${\cal W}^{-1}$ are related by:
$$W=\delta^u A_1A_2\cdots A_k, \ \
\ W^{-1} = \delta^{-(u+k)}\tau^{-(u+k)}(\ov{A_k})\cdots
\tau^{-(u+1)}(\ov{A_1}).$$
\end{enumerate}
\end{prop}

\proof
It is easy to show that the following identities hold for $A\in[0,1]$.
\begin{displaymath}
\begin{array}{ll}
S(\tau(A))=\tau(S(A)),  &R(\tau(A))=\tau(R(A)),\\
S(\ov{A})=R(A),     &R(\ov{A})=S(\tau(A)).
\end{array}
\end{displaymath}
Then (I) is clear because
$$R(\ov{B})\cap S(\tau(\ov{A}))=S(\tau(B))\cap
R(\tau(A)) =\tau(R(A)\cap S(B))$$
As for (II), it is easy to see that the equation for $W^{-1}$
holds. And it is the canonical form by (I). \endpf

\

We end this section with two technical lemmas and a corollary which will
play a  role in the implementation of Theorem \ref{thm:l canonical form} as
an algorithm.  They relate to the steps to be followed in the passage
from an arbitrary representative of a braid of the form
$\delta^uA_1A_2\cdots A_r$, where each $A_i$ is in $[0,1]$, to one in
which every adjacent pair $A_iA_{i+1}$ satisfies the conditions for
left-weightedness. The question we address is this: suppose that
$A_1A_2\dots A_i$ is left-weighted, that $A_{i+1}$ is a new canonical
factor, and that $A_iA_{i+1}$ is not left-weighted. Change to
left-weighted form $A_i^\prime A_{i+1}^\prime$, but now
$A_{i-1}A_i^\prime$ may not be left-weighted. We change it to
left-weighted form
$A_{i-1}^\prime A_i^{\prime\prime}$. The question which we address is
whether it is possible that after both changes
$A_i^{\prime\prime}A_{i+1}^\prime$ is not left-weighted?
The next two lemmas will be used to show that the answer is ``no".

\begin{lemm}
\label{lemm:AB,BC}
Let $AB, BC$ be canonical factors. Then
$A\lceil C$ if and only if $(AB)\lceil C$.
\end{lemm}

\proof By Corollary \ref{coro:lwtest},
$(AB)\lceil C$ iff for each $c\in S(C)$ there exists
$a\le (AB)$ such that
$(a,c)$ is an obstructing pair. Since $BC$ is a canonical factor,
we know from Corollary \ref{coro:BaCbD} that we cannot have $a\le B$.
Therefore the only possibility is that $a\le A$.
\endpf

\begin{lemm}\label{lemm:ABCD}
Suppose that $A,B,C,D,B',C'$ are canonical factors and that $ABCD=AC'B'D$.
Suppose also that $BC,CD,AC',C'B'$ are canonical factors, and that
$A\lceil B$ and $B\lceil D$. Then $B'\lceil D$.
\end{lemm}

\proof
By Proposition
\ref{prop:invform} it suffices to show that $\ov{D}\lceil
\tau(\ov{(B')})$.
\begin{eqnarray*}
S(\ov{D})&\subseteq & S(\ov{D}\tau(\ov{(B')}))\\
&\subseteq & S(\ov{D}\tau(\ov{(B')})\tau^2(\ov{(AC')}))\\
&=& S(((AC')B'D)^{-1}\delta^3)\\
&=& S((AB(CD))^{-1}\delta^3)\\
&=& S(\ov{(CD)}\tau(\ov{B})\tau^2(\ov{A}))\\
&=& S(\ov{(CD)}\tau(\ov{B})) \\
&=& S(\ov{D}\tau(\ov{(BC)}))\\
&=& S(\ov{D})
\end{eqnarray*}

\noindent Here the fourth equality which follows the first two
inclusions is a consequence of the fact that
$A\lceil B$. The sixth equality follows from $B\lceil D$, which (by Lemma
\ref{lemm:AB,BC}) implies that $(BC)\lceil D$.  But then, every inclusion
must be an equality, so that $S(\ov{D})=S(\ov{D}\tau(\ov{B'}))$. But then,
by Proposition  \ref{prop:invform}, it follows that $B'\lceil D$.
\endpf

We now apply the two lemmas to prove what we will need about
left-weightedness.

\begin{coro}
\label{coro:left-weightedness}
Suppose that $A_{i-1},A_i,A_{i+1}$ are canonical factors, with
$A_{i-1}\lceil A_i$.
Let $A_i^\prime, A_{i+1}^\prime$ be canonical factors with
$A_iA_{i+1} = A_i^\prime A_{i+1}^\prime$ and
$A_i^\prime\lceil A_{i+1}^\prime$.
Let $A_{i-1}^\prime,A_i^{\prime\prime}$ be canonical factors with
$A_{i-1}A_i^\prime = A_{i-1}^\prime A_i^{\prime\prime}$ and
$A_{i-1}^\prime\lceil  A_i^{\prime\prime}$. Then
$A_i^{\prime\prime}\lceil A_{i+1}^\prime$.
\end{coro}

\proof  The conversion of $(A_i)(A_{i+1})$ to left-weighted form
$A_i^\prime\lceil A_{i+1}^\prime$ implies the existence of $U\geq e$
with
$$(A_i)(A_{i+1}) = (A_i)(UA_{i+1}^\prime) = (A_iU)(A_{i+1}^\prime) =
(A_i^\prime)(A_{i+1}^\prime).$$ The
subsequent conversion of $(A_{i-1})(A_i^\prime)$ to left-weighted form
$A_{i-1}^\prime\lceil A_i^{\prime\prime}$ implies the existence of
$V\geq e$
with
$$(A_{i-1})(A_i^\prime) = (A_{i-1})(VA_i^{\prime\prime}) =
(A_{i-1}V)(A_i^{\prime\prime})=
(A_{i-1}^\prime)(A_i^{\prime\prime}).$$
Set $A_{i-1}=A, \ A_i=B, \ U=C, \ A_{i+1}^\prime = D,
\ V=B^\prime, \ A_i^{\prime\prime}=C^\prime$ and apply
Lemma \ref{lemm:ABCD}
to conclude that $A_i^{\prime\prime}\lceil A_{i+1}^\prime$.  $\|$
\newpage

\section{Algorithm for the word problem, and its complexity}
In this section we describe our algorithm for putting an arbitrary
${\cal W}\in B_n$  into left-canonical form and analyze the complexity of
each step in the algorithm. The {\it complexity} of a computation is said to
be ${\cal O}(f(n))$ if the number of steps taken by a Turing machine (TM) to
do the computation is at most $kf(n)$ for some positive real nymber $k$.
Our calculations will be based upon the use of a random access
memory machine (RAM), which is in general
faster than a TM model (see Chapter 1 of \cite{AHU}). An RAM machine has
two models: in the first (which we use) a single input (which we interpret
to be the braid index) takes one memory unit of time. Unless the integer
$n$ is so large that it cannot be described by a single computer word, this
`uniform cost criterion' applies. We assume that to be the case, i.e. that
the braid index
$n$ can be stored by one memory unit of the machine.

We recall that each canonical factor decomposes into a product of parallel
descending cycles $A=\delta_{\pi_1}\cdots\delta_{\pi_k},$ and that $A$ is
uniquely determined by the permutation
$\pi_1\cdots\pi_k$.
So we identify a canonical factor with the permutation
of its image under the projection $B_n\to\Sigma_n$.
We denote each cycle $\pi_i$ by its ordered sequence of subscripts.
For example, we write (5,4,3,1) for $a_{54}a_{43}a_{31}$.

We use two different ways to denote a permutation $\pi$ which is the image
of a canonical factor: the first is by the $n$-tuple
$(\pi(1),\ldots,\pi(n))$ and the second is by its decomposition as a
product of parallel, descending cycles
$\pi_1\cdots\pi_k$.
The two notations can be transformed to one another in linear time. The
advantage of the notation $\pi = (\pi(1),\ldots,\pi(n))$ is that the group
operations of multiplication and inversion can be perfomed in linear time.

If $A,B\in[0,1]$, the {\em meet} of $A$ and $B$,
denoted $A\wedge B$, is defined to be the maximal
canonical factor $C$ such that $C\le A$ and $C\le B$.
Our definition is analogous to that in \cite[page 185]{Epstein}.
Note that $C$ can be characterized by the property that $S(C) = S(A)\cap
S(B)$.

\begin{lemm}\label{thm:meet}
If $A,B\in[0,1]$, then $A\wedge B$ can be
computed in linear time as a function of $n$.
\end{lemm}

\proof
Let $A=\pi_1\cdots\pi_k$ and $B=\tau_1\cdots\tau_\ell$, where the ordering
of the factors is arbitrary, but once we have made the choice we
shall regard it as fixed. Let $\coprod$ denote disjoint union. Then
$A\wedge B=\prod_{i,j}\pi_i\wedge
\tau_j$ since
\begin{eqnarray*}
S(A\wedge B) &=& S(A)\cap S(B) = \bigg(\coprod_i S(\pi_i)\bigg)
     \cap\bigg( \coprod_jS(\tau_j) \bigg)\\
    &=& \coprod_{i,j} \left( S(\pi_i)\cap S(\tau_j)\right)
        =\coprod_{i,j} \left( S(\pi_i\wedge\tau_j)\right).
\end{eqnarray*}
For two descending cycles $\pi_i=(t_1,\ldots,t_p)$ and
$\tau_j=(s_1,\ldots,s_q)$, we have:
$$
S(\pi_i)\cap S(\tau_j)
 = \{a_{ts}\mid \mbox{$t>s$ and $t,s\in \{t_1,\ldots,t_p
\}\cap\{s_1,\ldots,s_q
\}$}\}
$$
Thus
$$\pi_i\wedge\tau_j=\left\{\begin{array}{ll}
(u_1,\ldots,u_r)\quad &\mbox{if }\{t_1,\ldots,t_p \}\cap
\{s_1,\ldots,s_q \}=\{u_1,\ldots,u_r\}
\mbox{ where } r\ge 2\\
e \quad &\mbox{if }|\{t_1,\ldots,t_p \}\cap\{s_1,\ldots,s_q \}|\le 1
\end{array}\right. $$
If we treat a decreasing cycle as a subset of $\{1,\ldots, n\}$ and a
canonical factor
as a disjoint union of the corresponding subsets, we may write
$A\wedge B$ as
$$A\wedge B=\coprod_{i,j} ( \pi_i\cap\tau_j).$$
We will find this disjoint union $A\wedge B$ of subsets of $\{1,\ldots,
n\}$ in linear time by the following four steps:

\begin{enumerate}
\item
Make a list of triples $\{(i,j,m)\}$ such that $m=1,\ldots,n$
appears in $\pi_i$ and $\tau_j$.
We do this by scanning $A=\pi_1\cdots\pi_k$ first and
writing $(i,\quad, m)$ if $\pi_i$ contains
$m$ and then scanning $B=\tau_1\cdots\tau_\ell$ and
filling in the middle entry of the triple
$(i,\quad, m)$ with $j$ if $\tau_j$ contains $m$.
We throw away all triples with a missing entry.
The list contains at most $n$ triples.
For example, if $A=(5, 4, 1)(3, 2)$ and $B=(4, 2, 1)$, our list
contains three triples, $(1,1,4)$, $(1,1,1)$, $(2,1,2)$.
This operation is clearly in ${\cal O}(n)$.

\item Sort the list of triples lexicographically.
In the above example, $(2,1,2)$, $(1,1,4)$, $(1,1,1)$ are the entries
in the sorted list.
There is an algorithm to do this in time ${\cal O}(n)$.
See \cite[Theorem 3.1]{AHU}.
\item Partition the sorted list by collecting triples with the same first
two entries and then throw away any collection with less than one element.
In the above example, $\{(2,1,2)\}$, $\{(1,1,4),(1,1,1)\}$ forms the
partitioned list and
we need to throw away the collection $\{(2,1,2)\}$.
This can be done by scanning the sorted list once.
Its complexity is ${\cal O}(n)$.

\item
>From each collection, write down the third entry to form a descending cycle.
Note that the third entries are already in the descending order.
In the above example, $(4,1)$ becomes $A\wedge B$. This step again
takes ${\cal O}(n)$.
\end{enumerate}
Since the above steps are all in ${\cal O}(n)$, we are done.
\endpf

\begin{remark}
\label{complexity} {\rm We remark that the key step in
both our computation and that in \cite{Epstein} is in the computation of
$A\wedge B$, where
$A$ and $B$ are permutation braids in \cite{Epstein} and canonical factors
in our work. Our computation is described in Lemma
\ref{thm:meet}. We now examine theirs. The set
$R_\sigma$ which is used in
\cite{Epstein} is defined on pages 184-5 of \cite{Epstein} and
characterized inductively on page 185. The fact that
$R_\sigma$ is defined inductively means that one cannot use a standard
merge-sort algorithm. To get around this, the
merge-sort approach is modified, as explained on lines 2-5 of page 206:
one sorts
$1,2,\dots,n$ by the rule $i < j$ if $C(i) < C(j)$, where $C(i)$ is the
image of $i$ under the permutation $C=A\wedge B$. This ordering
between two integers cannot be done in constant time.
The running time is
${\cal O}(n\log n)$, where $\log n$ is the depth of recursion and
$n$ is the time needed to assign integers as above and
to merge sets at each depth.}
\end{remark}

\begin{lemm}\label{thm:leftweightedalg}
Let $A$, $B$ be canonical factors, i.e. $e\le A, B\le\delta$.
There is an algorithm
of complexity ${\cal O}(n)$ that converts $AB$ into the left weighted
decompostion, i.e. $A\lceil B$.
\end{lemm}

\proof
Let $A^\ast$ be the right complement of $A$, i.e.
$A^\ast\in[0,1]$ and $AA^\ast = \delta$.
Let $C=A^\ast \wedge B$ and $B=CB'$ for some $B'\in[0,1]$.
Then $(AC)\lceil B'$, for if there is $a_{st}\le B'$ such that
$ACa_{st}\in[0,1]$, then
$a_{st}\le B'\le B$ and $a_{st}\le (AC)^\ast \le A^\ast$, which is impossible
by the definition of meet.
Thus the algorithm to obtain the left weighted decomposition
consists in the following four steps:
\begin{enumerate}
\item[(I)] Compute the right complement $A^\ast$ of $A$.
\item[(II)] Compute $C=A^\ast\wedge B$.
\item[(III)] Compute $B'$ such that $B=CB'$.
\item[(IV)] Compute $AC$.
\end{enumerate}

The step (II) is in ${\cal O}(n)$ by the lemma \ref{thm:meet}.
The steps  (I), (III) and (IV) are in ${\cal O}(n)$ since they
involve inversions and
multiplications of permutations like $A^\ast = A^{-1}\delta$ and
$B'=C^{-1}B$.
\endpf


Now the algorithm for the left canonical decomposition of
arbitrary words is given by the following four processes.

\

\noindent {\bf The algorithm.} We are given an element
${\cal W} \in B_n$ and a word $W$ in
the band generators which represents it.

\begin{enumerate}
\item
If $W$ is not a positive word, then the first step is to
eliminate each generator
which has a negative exponent, replacing it with
$\delta^{-1}A$ for some positive word $A\in[0,1]$.
The replacement formulas for the negative letters in $W$ is:
\begin{displaymath}
a_{ts}^{-1}=
\delta^{-1}(n,n-1,\dots,t+1,t,s-1,s-2,\dots,2,1)(t-1,t-2,\dots,s+1,
s)
\end{displaymath}
The complexity of this substitution process is at most ${\cal O}(n|W|)$.
Notice that $|P|$ can be as long as ${\cal O}(n|W|)$, because each
time we eliminate a negative letter we replace it by a canonical factor
of length $n-2$.

\item
Use the formulas:
\begin{displaymath}
A_i\delta^k = \delta^k \tau^k(A_i)\quad\mbox{and}\quad
\delta^{-1}\delta^k = \delta^{k-1}
\end{displaymath}
to move $\delta^{-1}$'s to the extreme left,
to achieve a  representative of $\cal W$ of the form
\begin{equation}\label{eqn:decomposition}
W = \delta^uA_1A_2\cdots A_k,\quad A_i\in[0,1]
\end{equation}
and $|A_i| = 1$ or $n-2$, according as $A_i$ came from a positive or a
negative letter in $W$.
Since we can do this process by scanning the word just once,
the complexity of this rewriting process depends
on the length of $A_1A_2\cdots A_k$ and so it is at most ${\cal O}(n|W|)$.

\item
Now we need to change the above decomposition
(\ref{eqn:decomposition}) to left
canonical form. In the process we will find that $u$ is maximized,
$k$ is minimized and $A_i\lceil A_{i+1}$ for every $i$ with
$1\le i\le k$. This can be achieved by repeated uses of
the subroutine that is described in the proof of Lemma
\ref{thm:leftweightedalg}.

In order to make the part $A_1A_2\cdots A_k$ left-weighted, we may
work either forward or backward.
Assume inductively that $A_1A_2\cdots A_i$ is already in its left
canonical form.
Apply the subroutine on $A_iA_{i+1}$ to achieve $A_i\lceil A_{i+1}$
and then to $A_{i-1}A_i$
to achieve $A_{i-1}\lceil A_{i}$. Corollary \ref{coro:left-weightedness}
guarantees that we still have
$A_i\lceil A_{i+1}$, i.e. we do not need to go back to
maintain the left-weightedness.
In this manner we apply the subroutine at most $i$-times
to make $A_1A_2\cdots A_iA_{i+1}$
left-weighted. Thus we need at most $k(k+1)/2$ applications of the
subroutine to complete the
left canonical form of $A_1A_2\cdots A_k$ and the complexity is ${\cal
O}(|W|^ 2n)$ since $k$
is proportional to $|W|$.

We may also work backward to obtain the same left canonical
form by assuming inductively that $A_iA_{i+1}\cdots A_k$ is already in its
canonical form and trying to make $A_{i-1}A_iA_{i+1}\cdots A_k$
left-weighted.

\item[4.] Some of canonical factors at the beginning of
$A_1A_2\cdots A_k$ can be $\delta$ and some
of canonical factors at the end of $A_1A_2\cdots A_k$ can be $e$.
These should be absorbed in
the power of $\delta$ or deleted.
Note that a canonical factor $A$ is $\delta$ if and only if $|A|=n-1$ and
$A$ is $e$ if and only if
$|A|=0$. Thus we can decide whether $A$ is $\delta$ of $e$ in ${\cal O}(n)$
and so the complexity of this process is at most
${\cal O}(kn)={\cal O}(|W|n)$.
\end{enumerate}

\begin{thm}
\label{thm:complexityword}
There is an algorithmic solution to the word problem that is ${\cal
O}(|W|^2n)
$ where
$|W|$ is the length of the longer word among two words in $B_n$ that are
being
 compared.
\end{thm}
\proof When we put two given words into their canonical forms,
each step has complexity at most ${\cal O}(|W|^2n)$.
\endpf

\newpage
\section{The conjugacy problem}
Let $W=\delta^uA_1A_2\cdots A_k$, be the
left-canonical  form of $W\in B_n$. The result
of a {\em cycling} (resp. {\em decycling}) of $W=\delta^uA_1A_2\cdots
A_k$, denoted by $\cy(W)$ (resp. $\dy(W)$), is the braid
$\delta^uA_2\cdots A_k\tau^{-u}(A_1)$ (resp. $\delta^u\tau^{u}(A_k)A_1\cdots
A_{k-1}$).  Iterated cyclings are defined recursively by $\cy^i(W) =
\cy(\cy^{i-1}(W),$ and similarly for iterated decyclings.
It is easy to see that both cycling and decycling do not decrease(resp.
increase) the inf(resp. sup).

With essentially no new work, we are able to show that the solution
to the conjugacy problem of \cite{Garside} and \cite{Elrifai:Morton}
can be adapted to our new presentation of $B_n$. This approach was
taken in
\cite{KKL} for $n=4$. But there are no new difficulties encountered when
one goes to arbitrary $n$. The following two theorems are the keys to the
solution to the conjugacy problem.

\begin{thm} {\rm (\cite{Elrifai:Morton},\cite{KKL})}
\label{thm:finding inf and sup}
Suppose that $W$ is conjugate to $V$.
\begin{enumerate}
\item[\rm (I)] If \ $\inf (V)>\inf (W)$,
then repeated cyclings will produce $\cy^l(W)$ with \\
$\inf (\cy^l(W))>\inf (W)$.
\item[\rm (II)] If \ $\sup (V)<\sup (W)$, then repeated decyclings will
produce
$\dy^l(W)$ with \\  $\sup (\dy^l(W))<\sup (W)$.
\item[\rm (III)]  In every conjugacy class, the maximum value of \ $\inf (W)$
and the
minimum value of \  $\sup (W)$ can be achieved simultaneously.
\end{enumerate}
\end{thm}

Theorem \ref{thm:finding inf and sup} tells how to find $\inf(V)$ and
$\sup(V)$, and a special set of words which are conjugate to $V$ and have
maximal $inf$ and minimal $sup$.
The next theorem tells how to find {\em all} words which
are conjugate to the given word and have those values of
$\inf$ and $\sup$:

\begin{thm}{\rm (\cite{Elrifai:Morton},\cite{KKL})}
\label{thm:seq}
Suppose that two $n$-braids $V,W\in [u,v]$ are in the same conjugacy
class. Then there is a sequence of $n$-braids
$V=V_0,V_1,\ldots,V_k=W$, all in $[u,v]$, such that each
$V_{i+1}$ is the conjugate of $V_i$ by some element of $[0,1]$.
\end{thm}

\

\noindent {\bf An algorithm for the solution to the conjugacy problem:}
\ We can now describe our solution to the conjugacy problem.  Suppose  that
two words $V,W$ represent conjugate elements ${\cal V},{\cal W}$  of $B_n$.
Recall the definitions of $inf(V)$ and $sup(V)$ which were given after
Proposition \ref{prop;partial order}.  By Theorem
\ref{thm:finding inf and sup}, $\inf(V)\le\sup(W)$ and
$\inf(W)\le\sup(V)$.   Let $u=\min\{\inf(V),\inf(W)\}$ and let
$v=\max\{\sup(V),\sup(W)\}$. Then $V,W\in[u,v]$. The canonical lengths $\sup
(V)-\inf(V)$ and $\sup(W)-\inf(W)$ are proportional to the word lengths
$|V|$ and $|W|$, respectively. Thus $v-u$ is at most ${\cal O}(|V|+|W|)$.
The cardinality $|[u,v]|$ is given by $|[0,1]|^{v-u}$.
Since $|[0,1]|\le 4^n$, it follows that $|[u,v]|$ is at most ${\cal
O}(\exp(n(|V|+|W|)))$. By Theorem \ref{thm:seq}, there is a sequence
$V=V_0,V_1,\ldots,V_k=W$ of words in $[u,v]$ such that each  element is
conjugate to the next one by an element of $[0,1]$. The length $k$ of this
sequence can be
$|[u,v]|$ (in the worst case) and so we must have $U^{-1}VU=W$ for some
positive word $U$ of canonical length $\le |[u,v]|$. Since there are
$|[0,1]|^{|[u,v]|}$ many positive words of canonical length $\le |[u,v]|$,
the number of all possible $U$ is at most
${\cal O}(\exp(\exp(n(|V|+|W|))))$.
This certainly gives a finite algorithm for the conjugacy problem.

\

A more sensible approach is as follows:
Given a $n$-braid $W$, the collection of conjugates of $W$ that has
both the maximal infimum and the minimal supremum is called
the {\it super summit set} of $W$ after
\cite{Garside},\cite{Elrifai:Morton}. Clearly the super summit set of a word
is an invariant of its conjugacy class. If we iterate the cycling operation
on a word $W$, then the fact that the number of positive words of fixed
length is finite insures that we eventually obtain positive integers $N,K$
such that $\cy^N(W)=\cy^{N+K}(W)$.  In view of Theorem \ref{thm:finding inf
and sup} we conclude that $\inf(\cy^N(W))$ is the maximum value of infimum
among all conjugates of $W$,  Similarly, by interated decycling on
$\cy^N(W)$, we have
$\dy^M\cy^N(W)=\dy^{M+L}\cy^{N+K}(W)$ and so we conclude that
$\sup(\dy^M\cy^N(W))$ is the minimum value of supremum among all all
conjugates of
$W$. Therefore $\dy^M\cy^N(W)$ belongs to the super summit set of $W$.
In order to decide whether two words $V,W$ in $B_n$ are conjugate, we
proceed as follows:
\begin{enumerate}
\item  Do iterated cycling and decycling on $V$ and $W$ until we have $V'$
and $W'$ in their super summit sets, respectively.
If $\inf(V')\neq \inf(W')$ or $\sup(V')\neq \sup(W')$ we
conclude that they are not conjuagate.

\item  If $\inf(V')=\inf(W')$ and $\sup(V')= \sup(W')$  they may still not
be conjugate. We must compute the entire super summit set of $V$ by using
Theorem \ref{thm:seq} and the finiteness of the super summit set.

\item If any one element in the super summit set of $W$, say $W'$, is also
in the super summit set of $V$, then $W$ and $V$ are
conjugate. Otherwise, they are not conjugate.
\end{enumerate}
In the worst case, this algorithm is no different from the previous one.
Nevertheless, we have lots of data which gives evidence
of additional structure, but we need to do more
work before we can improve the algorithm.

\

\noindent {\bf Example:} The conjugacy classes of the 4-braids which are
defined by the two words
$X,Y$ which are given below have the same `numerical class invariants', i.e.
the same $inf$,
$sup$ and cardinality of the super summit set. The super summit sets split
into orbits under cycling and decycling, and the numbers and lengths of
these orbits coincide. But the braids are not conjugate because their super
summit sets are disjoint:
$$X=a_{43}^{-2}a_{32}a_{43}^{-1}a_{32}a_{21}^3a_{32}^{-1}
a_{21}a_{32}^{-1}, \ \ \ \ \
Y=a_{43}^2a_{32}^{-1}a_{21}^3a_{32}a_{43}^{-1} a_{21}^{-1}a_{32}^{-2}.$$

\newpage

\noindent
\footnotesize
{Joan Birman, Address: Department of Mathematics, Columbia
University, New York, NY 10027.\quad Email: jb@math.columbia.edu,
Institutional affiliation: Barnard College of Columbia University. \\
\ \\
Ki Hyoung Ko, Address: Department of Mathematics, Korea Advanced
Institute of Science and Technology, Taejon, 305-701, Korea.\quad
Email: knot@knot.kaist.ac.kr\\
\ \\
Sang Jin Lee, Address: Department of Mathematics, Korea Advanced
Institute of Science and Technology, Taejon, 305-701, Korea.\quad
Email: sangjin@knot.kaist.ac.kr

\end{document}